%% file: ms.tex
\documentclass{amsart}

\usepackage[utf8]{inputenc}
\usepackage[T1]{fontenc}
\usepackage{amsmath}
\usepackage{amsfonts}
\usepackage{amssymb}
\usepackage{graphicx}
\usepackage{color}
\usepackage{ulem}
\usepackage{selinput}
\usepackage{caption}
\usepackage{subcaption}
\usepackage{cancel}
\usepackage{enumitem}
\usepackage[paperwidth=7in, paperheight=10in, margin=0.8in]{geometry}
\usepackage[backref,colorlinks,linkcolor=red,anchorcolor=green,citecolor=blue]{hyperref}
\usepackage{todonotes}
\usepackage{stmaryrd}
\usepackage{tikz}
\usepackage{booktabs}

\usetikzlibrary{calc}

\newtheorem{proposition}{Proposition}[section]

\theoremstyle{remark}
\newtheorem{remark}[proposition]{Remark}
\DeclareMathOperator{\Div}{div}

\author{Fischer Michael}
\address{Radon Institute for Applied and Computational Mathematics, Altenbergerstr. 69, 4040 Linz, Austria}
\email{tbc}
\author{Jankowiak Gaspard}
\address{Radon Institute for Applied and Computational Mathematics, Altenbergerstr. 69, 4040 Linz, Austria}
\author{Wolfram Marie-Therese}
\address{University of Warwick, Mathematics Institute, Gibbet Hill Road, CV47AL Coventry, UK and Radon Institute for Applied and Computational Mathematics, Altenbergerstr. 69, 4040 Linz, Austria}

\title[Crowding and pushing in corridors]{Micro- and macroscopic modeling of crowding and pushing in corridors}
\date{}

\newtheorem{theorem}{Theorem}

\begin{document}
\maketitle

\input{document}

\bibliographystyle{abbrv}
\bibliography{bibliography}
\end{document}

%% file: document.tex
%The abstract of your paper
\begin{abstract}
Experiments with pedestrians revealed that the geometry of the domain, as well as the incentive of pedestrians to reach a target as fast as possible
have a strong influence on the overall dynamics. In this paper, we propose and validate different mathematical models
at the micro- and macroscopic levels to study the influence of both effects. We calibrate the models with experimental data and compare the results at the
micro- as well as macroscopic levels. Our numerical simulations reproduce qualitative experimental features on both levels, and indicate
how geometry and motivation level influence the observed pedestrian density.
Furthermore, we discuss the dynamics of solutions for different modeling approaches and comment on the analysis of the respective equations.
\end{abstract}

\section{Introduction}

\noindent In this paper, we develop and analyze mathematical models for crowding and queuing at exits and bottlenecks, which are motivated by
experiments conducted at the Forschungszentrum J\"ulich and the University of Wuppertal, see~\cite{art:armin}. In these experiments, student groups
of different size were asked to exit through a door as fast as possible. Each run corresponded to different geometries of the domain, ranging from
a narrow corridor to an open space, as well as different motivation levels, by giving more or less motivating instructions.
The authors observed that
\begin{itemize}[leftmargin=10pt,topsep=5pt]
 \item The narrower the corridor, the more people lined up. This led to a significantly lower pedestrian density in front of the exit.
 \item A higher motivation level led to an increase of the observed densities. However its impact on the density was smaller than changing the shape
 of the domain.
\end{itemize}
Adrian et al.~\cite{art:armin} supported their results by a statistical analysis of the observed data as well as computational experiments using a force based
model. We follow a different modeling approach in this paper, proposing and analyzing a cellular automaton (CA) model which is motivated by the aforementioned experiments.
We see that these minimalistic mathematical models reproduce the observed behavior on the microscopic as well as macroscopic level.

 There is a rich literature on mathematical models for pedestrian dynamics. Ranging from microscopic agent or cellular automaton based
approaches to the macroscopic description using partial differential equations.
The social force model, see~\cite{HT77, O79, HM1995}, is the most prominent individual based model. Here pedestrians are characterized by their position and velocity,
which change due to interactions with others
and their environment. More recently, the corresponding damped formulation, see~\cite{art:armin}, has been considered in the literature.
In cellular automata (CA), another much used approach, individuals move with given rates from one discrete cell to another. One advantage of
CA approaches is that the formal passage
from the microscopic to the macroscopic level is rather straight-forward based on a Taylor expansion of the respective transition rates. This
can for example be done systematically using tools from symbolic computation, see~\cite{KRRW2015}. CA approaches have been used successfully to describe
lane formation, as for example in~\cite{NS2012}, or evacuation situations, such as in~\cite{art:conflict3}. The dynamics of the respective
 macroscopic models was investigated in various situations such as uni- and bidirectional flows or cross sections, see for example~\cite{BHRW2015, art:burger}.

 Macroscopic models for pedestrian dynamics are usually based on conservation laws, in which the average velocity of the crowd is reduced due to interactions with
others, see~\cite{a:Piccoli2011, a:macro3}. In general it is assumed that the average speed changes with the average pedestrian density, a relation known as
the fundamental diagram. In this context, finite volume effects, which ensure that the maximum pedestrian density does not exceed a certain physical bound, play an
important role. These effects result in nonlinear diffusivities, which saturate as the pedestrian density reaches the maximum density,
and cross-diffusion in case of multiple species, see for example~\cite{BHRW2015}.
One of the most prominent macroscopic models is the Hughes model, see~\cite{a:hughes1, a:Hughes0}.
It consists of a nonlinear conservation law for the pedestrian density which is coupled to the eikonal equation to determine the shortest path to a target (weighted by the
pedestrian density). We refer to the textbooks by Cristiani et al., see~\cite{book:tosin} and Maury and Faure, see~\cite{book:maury}, for a more detailed overview on pedestrian dynamics.

 Many PDE models for pedestrian dynamics can be interpreted as formal gradient flows with respect to the Wasserstein distance. In this context, entropy methods
have been used successfully to analyze the dynamics of such equations. For example, the boundedness by entropy principle ensures the global in time existence of weak solutions for large classes of nonlinear partial differential equation systems, see~\cite{J2015}.
These methods have been proven to be useful also in the case of nonlinear boundary conditions and were also used by Burger and Pietschmann~\cite{art:burger} to show
existence of stationary solutions to a nonlinear PDE for unidirectional pedestrian flows with nonlinear inflow and outflow conditions. The respective
time dependent result was subsequentially presented in~\cite{GSW2019}.

 The calibration of microscopic pedestrian models is of particular interest in the engineering community. Different calibration techniques have been
used for the social force model, see~\cite{a:cal2, a:cal1} and CA approaches, see~\cite{a:cal7, a:cal6}. Nowadays a large amount of data is publicly available - for example
the database containing data for a multitude of experimental setups at the Forschungszentrum in J\"ulich, or data collected in a Dutch railway stations over the course
of one year, see~\cite{a:stat2}. However, many mathematical questions concerning the calibration of macroscopic and mean-field models from individual trajectories are still open.

 In this paper, we develop and analyze mathematical models to describe queuing individuals at exits and bottlenecks. Our main contributions are as follows:
\begin{itemize}[leftmargin=10pt,topsep=5pt]
 \item Develop microscopic and macroscopic models to describe pedestrian groups with different motivation levels and analyze their dynamics for various geometries.
 \item Calibrate and validate the microscopic model with experimental data in various situations.
 \item Compare the dynamics across scales using computational experiments.
 \item Present computational results, which reproduce the experimentally observed characteristic behavior.
\end{itemize}
\noindent This paper is organized as follows. We discuss the experimental setup and the proposed CA approach in Section~\ref{s:experiment_model}. In Section~\ref{s:valcal},
we present the details of the corresponding CA implementation and use experimental data to calibrate it. Section~\ref{s:pde} focuses on the description on the macroscopic level
by analyzing the solutions to the corresponding formally derived PDE. We conclude by discussing alternative modeling approaches in Section~\ref{s:alt} and summaries our
findings in Section~\ref{s:con}.

\section{The experimental setup and the microscopic model}\label{s:experiment_model}

\subsection{The experimental setup}

We start by discussing the experiments, which serve as the motivation for the proposed microscopic model, see~\cite{art:armin}. These experiments were
conducted at the University of Wuppertal, Germany. The respective data is available online, see~\cite{hp:juelich}.

\begin{figure}[htb]
 \centering
 \centering
\includegraphics[height=7cm]{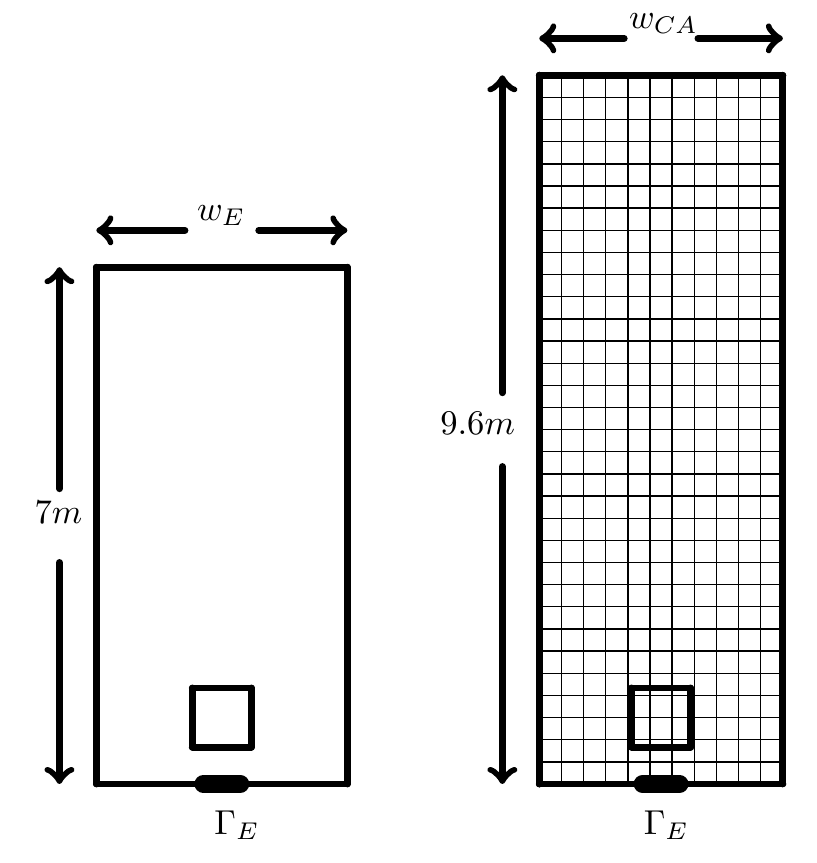}
\caption{Left: Sketch of experimental setup at the University of Wuppertal, showing the corridor width $w_E$ in the experiments,
the exit $\Gamma_E$ and the measurement area. Right: computational domain with adapted width $w_{CA}$ to ensure a consistent
discretization of the exit and an increased length $l_{CA}=9.6$m.}\label{fig:domaindisc}
\end{figure}

 The conducted experiments were designed to obtain a better understanding how social cues and the geometry of the domain influence individual
behavior. For this purpose runs with five different corridor widths, varying from $1.2$ to $5.6$ meters, were conducted over the course of several days. For each corridor,
a group of students was instructed to reach a target. These runs were then repeated with varying instructions, for example suggesting that queuing
is known to be more efficient or suggesting to go as fast as possible. The instructions were given to vary the motivation level and see their effect on the crowd dynamics.
The number of students in the different runs (which corresponded to the different corridor width) varied from
$20$ to $75$. The trajectory of each individual was recorded and used to compute the average density with the software
package JuPedSim, available at~\cite{hp:jupedsim}. The post-processed data showed that the average pedestrian density becomes particularly high in a $0.8\times0.8$ meter area,
$0.5$ meters in front of the exit, highlighted in Figure~\ref{fig:domaindisc}. Within this area, average densities up to $10$ $\text{p}/\text{m}^2$ (pedestrians per square meter) were observed. The densities varied significantly for the different
runs - they were much
lower for narrow, corridor-like domains and increased with the motivation level. Further details on the experimental setup can be found in ~\cite{art:armin}.

\subsection{The cellular automaton approach}\label{sec:ca}

In the following, we introduce a cellular automaton approach to describe the dynamics of agents queuing in front of the bottleneck. The dynamics of agents is determined by transition
rates, which depend on the individual motivation level and the distance to the target.

\begin{figure}[htb]
 \centering
\begin{subfigure}[b]{0.45\textwidth}
\centering
\includegraphics[height=5cm]{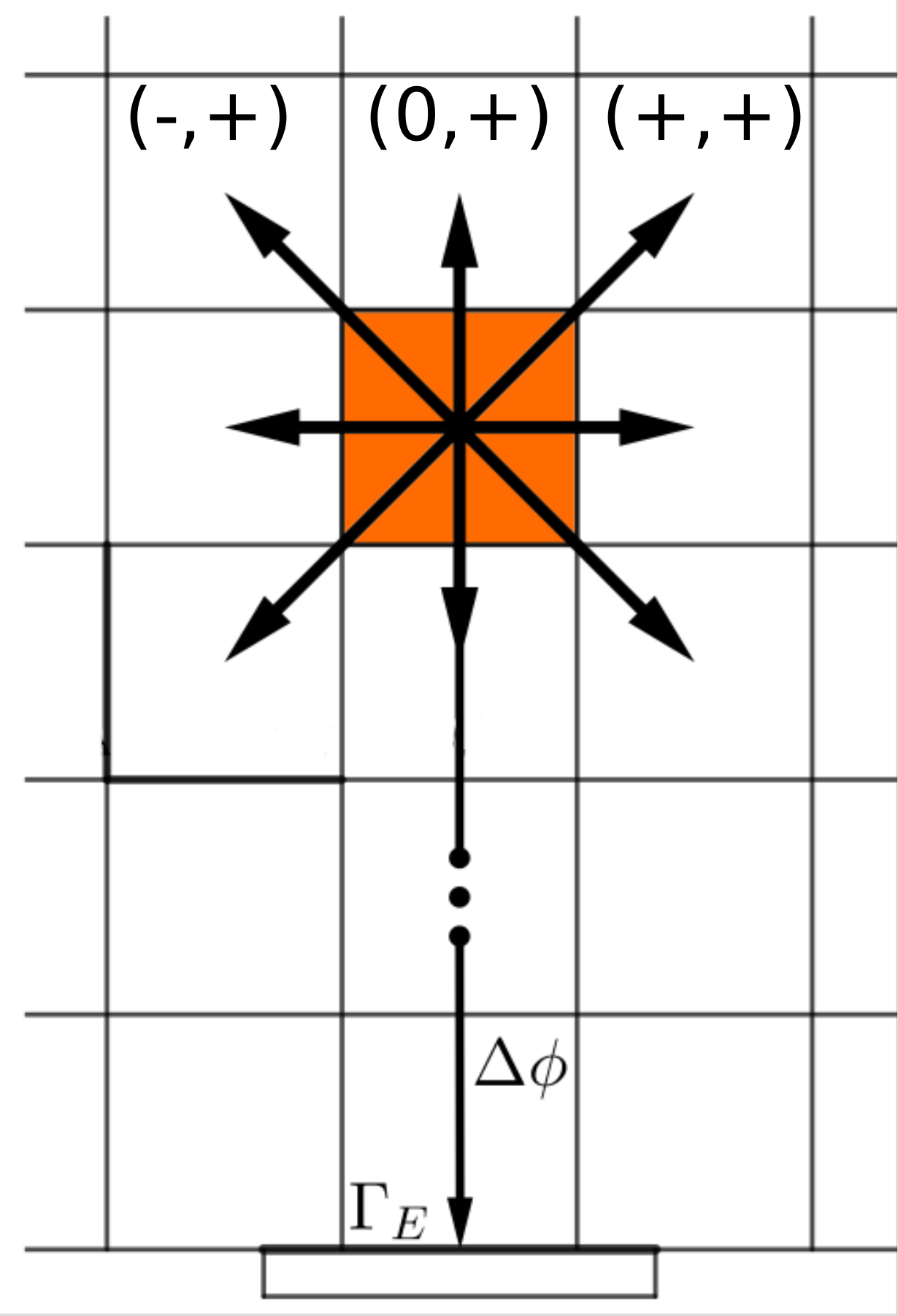}
 \caption{Moore neighborhood and selected elements of $I$; the potential $\phi$ corresponds to the distance to the exit.}\label{fig:moore}
\end{subfigure}
\hfill
\begin{subfigure}[b]{0.45\textwidth}
 \includegraphics[height=4cm]{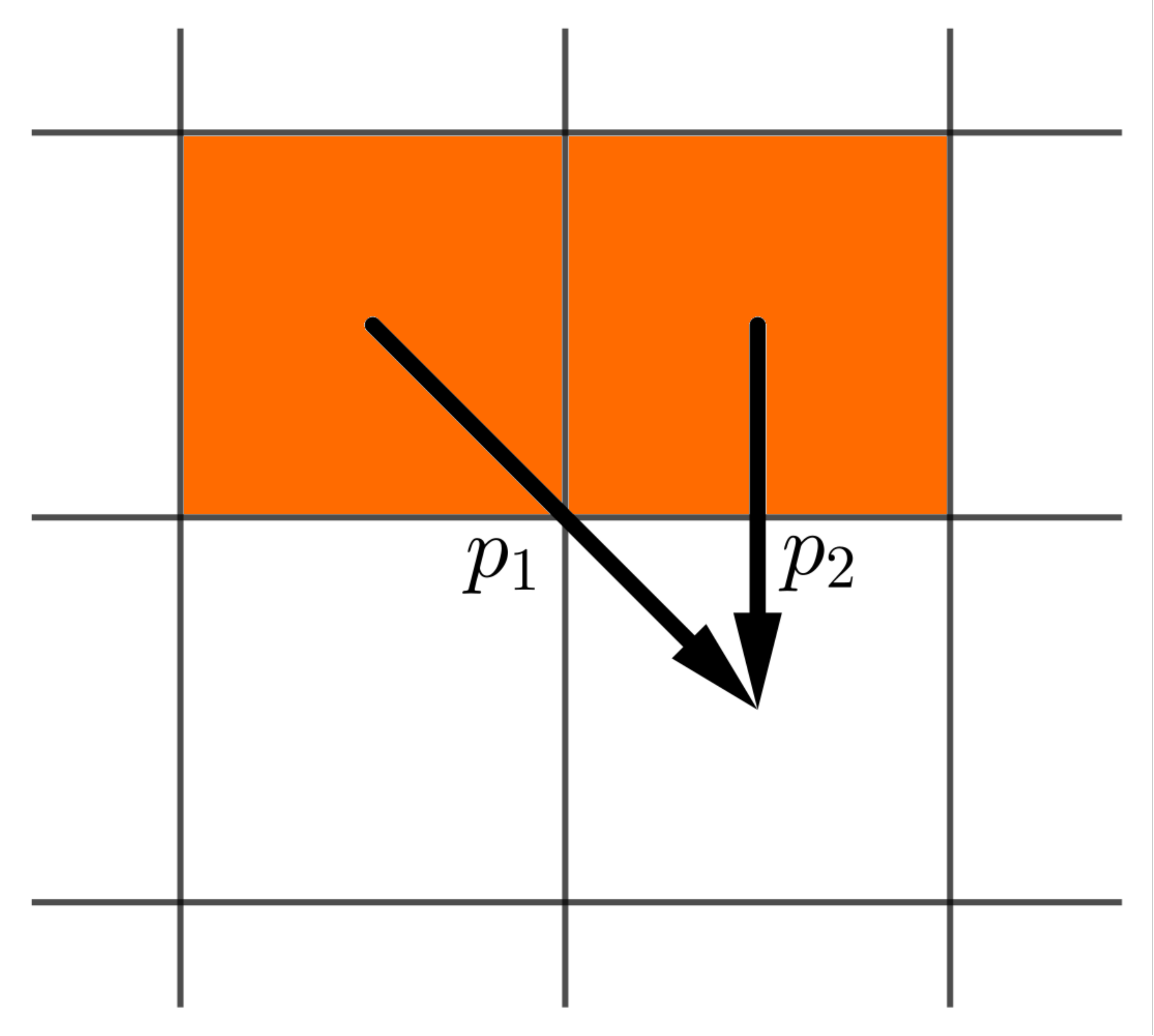}
 \caption{Conflict: Agent 1 wins with probability $\tilde p_1=p_1/(p_1+p_2)$, agent 2 with probability $1-\tilde p_1$.} \label{fig:conflict}

\end{subfigure}
\caption{Cellular automaton: transition rules.}
\end{figure}

We split the domain into squares with sides of length $\Delta x = 0.3$m. This discretization corresponds to a maximum
packing density of $11.11$ $\text{p}/\text{m}^2$. Cell-sizes of $0.09$m$^2$ have a comparable area to
ellipses with semi-axes $a=0.23$m and $b=0.12$m - a reference measure
for pedestrians commonly used in
agent based simulations, see~\cite{mt:ben}. The cellular automaton is implemented on a Moore neighborhood, see Figure~\ref{fig:moore}. Agents are allowed to move into the eight neighboring sites.
Their transition rates
depend on the availability of a site -- a site can only be occupied by a single agent at a time -- the potential $\phi$, which corresponds to the minimal distance to the exit, as well
as the individual motivation level.
The positions of all agents are updated simultaneously, which is known as a parallel update. To do so, we calculate the transition rates for every agent
and resolve possible conflicts. In case of a conflict, the respective probabilities of the two agents wanting to move into the site are re-weighted, and one of them is selected.
This solution has been proposed by ~\cite{art:conflict2, art:conflict3} and is illustrated in Figure~\ref{fig:conflict}.

\begin{figure}[htb]
 \centering
\begin{subfigure}[b]{0.39\textwidth}
\centering
\includegraphics[height=4cm]{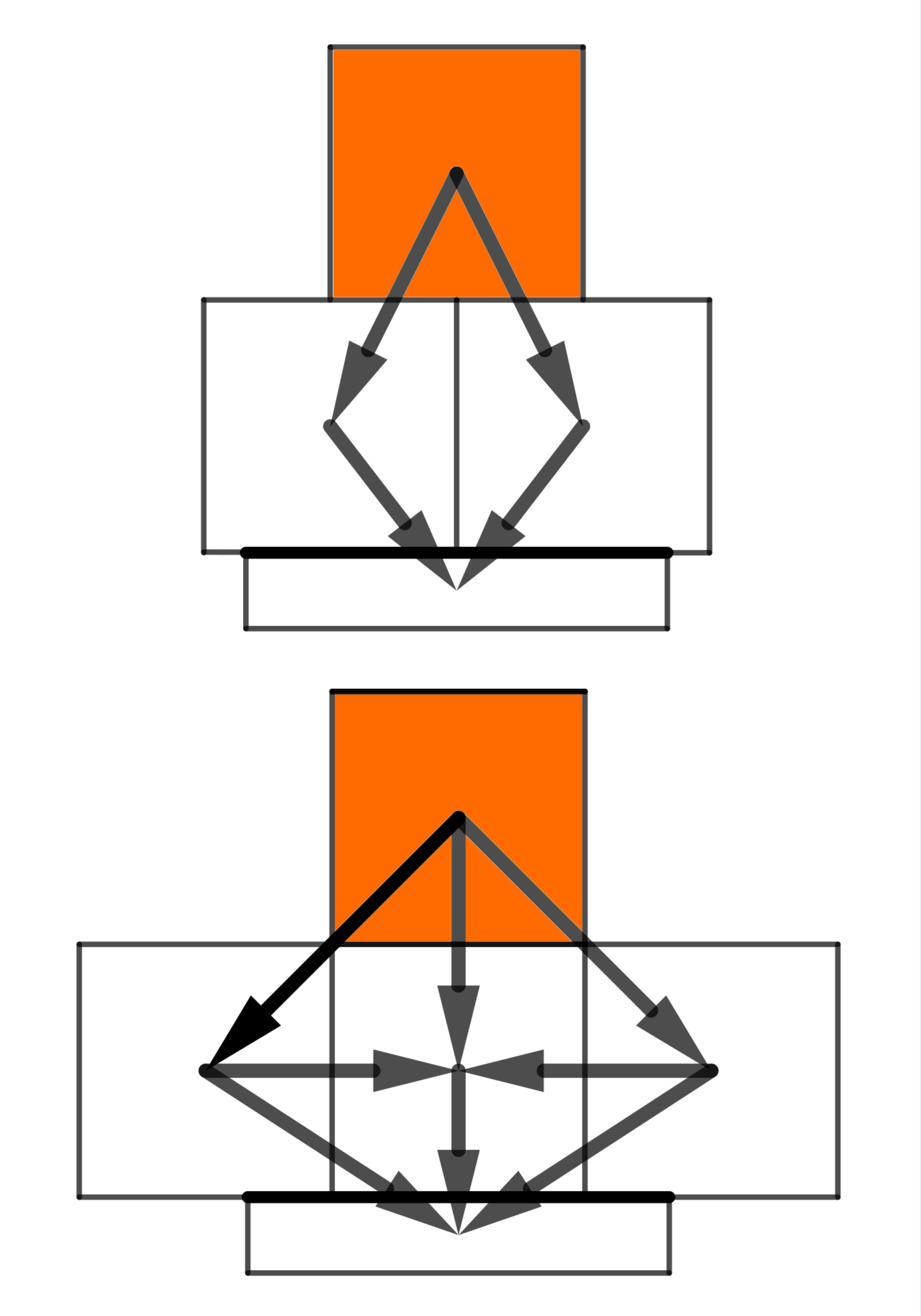}
 \caption{Two vs. three cells.}\label{fig:markov1}
\end{subfigure}
\begin{subfigure}[b]{0.6\textwidth}
 \centering
 \includegraphics[height=3.5cm]{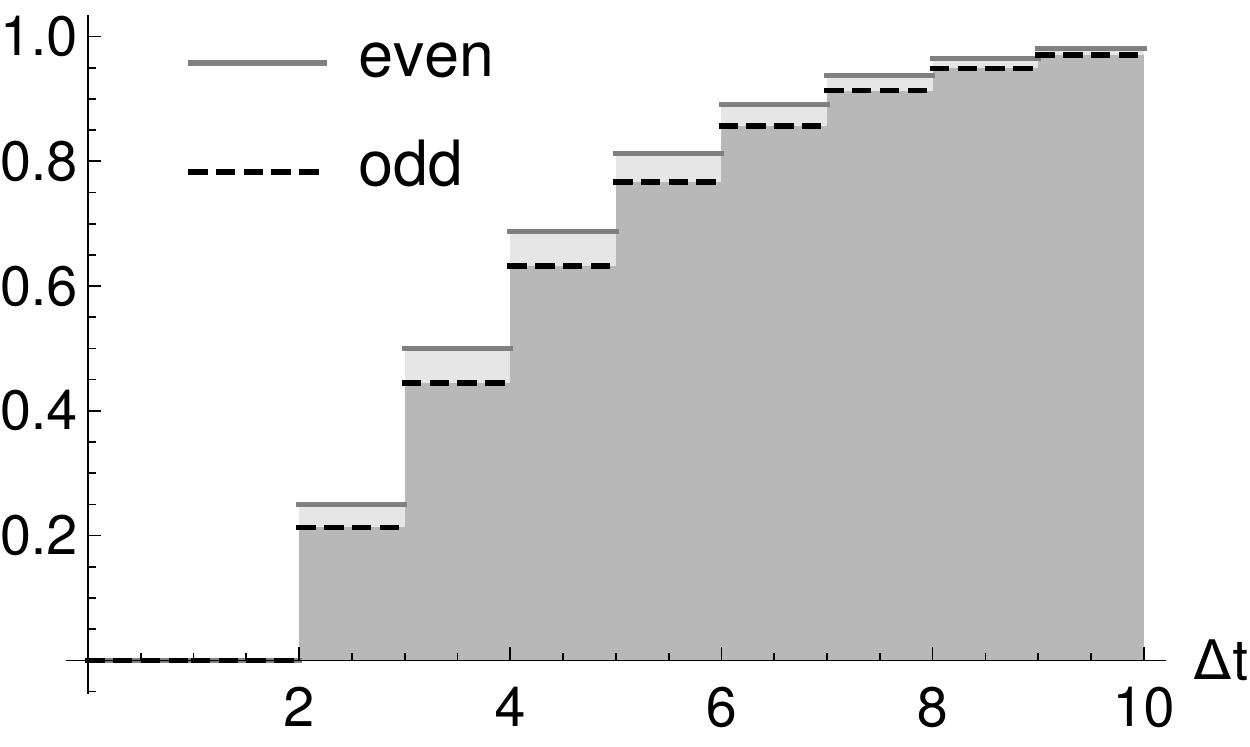}
 \caption{Two vs. three: comparison of the exit rate as a function of time.} \label{fig:markov2}
\end{subfigure}
\caption{Discretization of the exit. In the case of an even number of cells, a central positioned agent will leave the corridor faster than in the case of three cells. }
\end{figure}

Particular care has to be taken when modeling the exit. In doing so, we consider the special Markov-process, where a single agent is located at distance $\Delta x$ to the exit, see Figure~\ref{fig:markov1}.
We see that the exit can stretch over two or three cells. However, each setting has different exit probabilities and influences the exit rate. Figure~\ref{fig:markov2} illustrates the
different exit rates for the two situations in case of a single agent. We observe that the exit rate is higher if the exit is discretized using two cells.
To ensure a consistent discretization of the exit for all corridor widths, we choose a discretization using three cells for all corridors. Therefore, we changed the respective corridor
widths in the presented computational experiments from $1.2$m, $3.4$m and $5.6$m to $0.9$m, $3.3$m and $5.7$m, as illustrated in Figure~\ref{fig:domaindisc}. Furthermore, we extended the corridor to $9.6$m to ensure
sufficient space for all agents in case of
 larger groups. Note that in the actual experiments individuals were waiting behind the corridor entrance.

\subsubsection*{Transition-rates and the master equation}
The transition rates are based on the following assumptions, which are motivated by the previously detailed experiments:
\begin{itemize}[leftmargin=10pt,topsep=5pt]
 \item Individuals want to reach a target as fast as possible.
 \item They can only move into a neighboring site if it is not occupied.
 \item The higher the motivation level, the larger the transition rate.
\end{itemize}
Let $\rho = \rho(x,y,t)$ denote the probability of finding an individual at site $(x,y)$ at time $t$ and let $\mu$ denote the motivation level.
We will use the following abbreviation to state the master equation in 2D. Let $I$ denote the Moore neighbourhood of the cell $(x,y)$; then the
neighboring cells are indexed using the signs $I:=\{-,0,+\}^2\setminus \{(0,0)\}$, see Figure \ref{fig:moore}. The transition rate is given by
\begin{align}\label{eq:transitionrate}
\begin{split}
\mathcal T^{ij}(x,y)
&= \frac{1}{8(3-\mu)}
 \exp \left( \beta(\phi (x,y)-\phi (x+i\Delta x,y+j\Delta x))\right).
 \end{split}
\end{align}
The prefactor $\frac{1}{8}$ is a scaling constant such that
\begin{align*}
\sum_{(i,j)\in I\cup \{(0,0)\}} \mathcal T^{ij}(x,y)=1+\mathcal{O}((\Delta x)^2)
\end{align*}
holds. The parameter $\beta$ plays an important role to weigh the transition rates to the neighboring sites.
For $\beta = 0$, the transition rates are equidistributed over the neighboring cells and therefore the dynamics would
correspond to a random walk. In the limit $\beta~\rightarrow~\infty$, individuals will move in direction of the steepest
descent of $\phi$, and the dynamics become deterministic. Since the potential $\phi$ corresponds to the distance to the
target, the parameter $\beta$ has to scale as m$^{-1}$. Maury suggests in \cite{book:maury} that it should
be proportional to the characteristic distance  $\Delta x^{-1}$, which would correspond to the value $3.33$ in our setting.
The prefactor $(3-\mu)^{-1}$ changes the transition rates depending on the motivation level $\mu \in (-\infty, 1)$.
The smaller $\mu$, the less likely an agent is to move, see Remark \ref{rem:jumping}.
Additionally we consider size-exclusion, which corresponds to the prefactor $(1-\rho (x+i\Delta x,y + j \Delta x,t))$ in the following master-equation. It
ensures that the target site is not already occupied. Then the probability that a site $(x,y)$ is occupied at time $t + \Delta t$, is given by
\begin{multline}
    \label{eq:master}
\rho(x,y,t+\Delta t)= \rho(t, x,y)\\-\rho (x,y,t)
\sum_{(i,j)\in I} \mathcal T^{ij}(x,y) (1-\rho (x+i\Delta x,y+j\Delta x,t)) \\
 + \sum_{(i,j)\in I}\rho (x+i\Delta x,y+j\Delta x,t) \mathcal T^{ij}(x+i\Delta x,y+j\Delta x) (1-\rho (x,y)).
\end{multline}
In short, the first sum corresponds to all possible moves of an agent in $(x,y)$ to neighboring sites. The second sum all possible moves
from neighboring agents into that site.

We recall that agents can leave the domain from all three fields in front of the exit. In a possible conflict situation, that is two or three agents
located in the exit cells want to leave simultaneously, the conflict situation is resolved and the winner exists with
probability $p_{ex}$.

\begin{remark}
The choice of a Moore neighborhood instead of a Neumann neighborhood (as in~\cite{ NS2012,a:cavel}), is based on the experimental observations (individuals make diagonal
moves to get closer to the target). However, the choice of the neighborhood does not change the structure of the limiting partial differential
equation.
\end{remark}

\begin{remark}\label{rem:jumping}
Note that for the largest motivation level, that is $\mu=1$ the probability of staying is given by
\begin{align*}
\mathcal{T}^{00}(x,y)=1-\sum_{(i,j)\in I}\mathcal T^{ij}(x,y)=\frac{2-\mu}{3-\mu}=\frac{1}{2}.
\end{align*}
Such agents will move every second time-step. We see that the motivation $\mu$ has a direct influence on the desired maximum velocity $v_{max}$ on a microscopic level.
It also ensures that it is very unlikely that individuals step back in the case of a high number of agents between the agent and its target $\Gamma_E$.
\end{remark}

\section{Validation and calibration of the CA model}\label{s:valcal}

\subsection{Implementation of the CA approach}
We start by briefly discussing the implementation of the CA, which will be used for the calibration in the subsequent section. A CA simulation returns the average exit time
(that is the time when the last agent leaves the corridor) depending on the number of agents $n$, the corridor-width $w\in \{ 0.9, 3.3, 5.7\}$,
the motivation $\mu$, the length of a time-step $\Delta t$ and the parameter $\beta$. Each CA simulation is initialized with a random uniform distribution of agents.
For given parameters the returned average exit time $\bar{T}$ and maximum observed density is estimated by averaging over $5000$ CA simulations. Note that we calculate
this density in the area highlighted in Figure~\ref{fig:domaindisc}.

\begin{figure}[htb]
 \centering
\begin{subfigure}[t]{.19\textwidth}
 \centering
 \includegraphics[height=5cm]{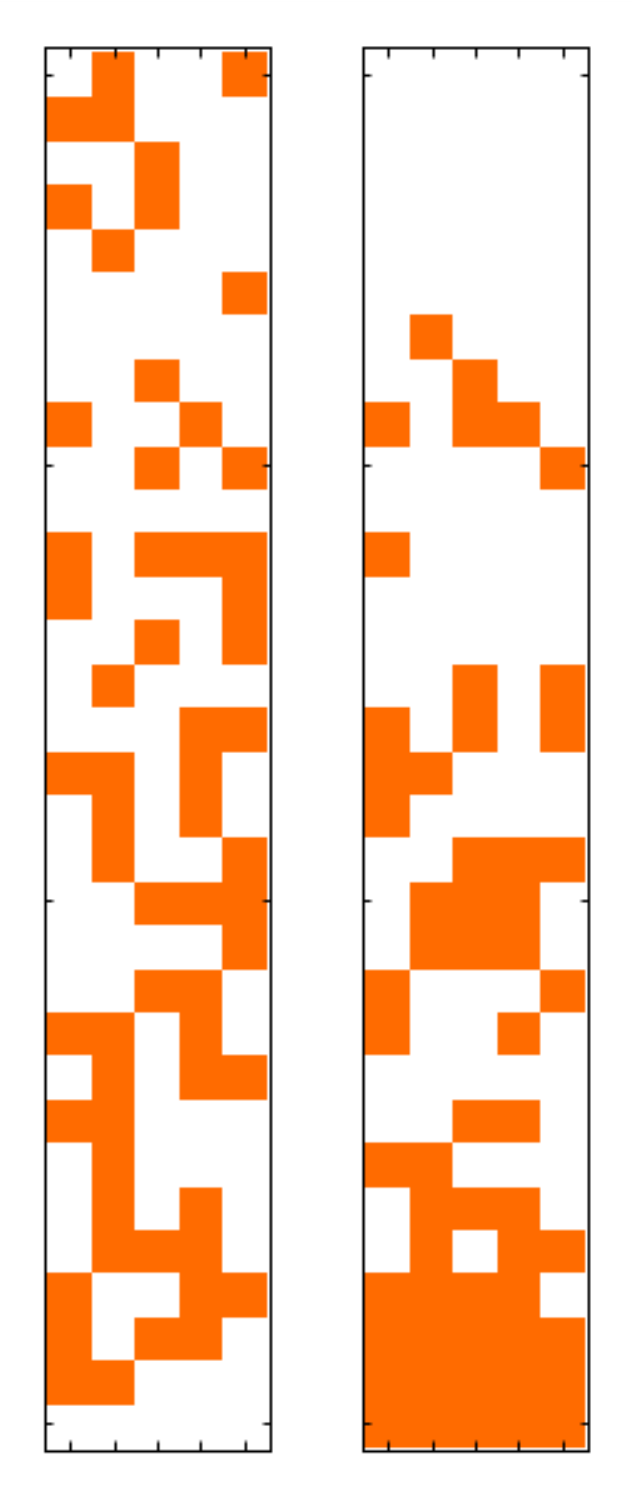}
 \caption{}
 \label{fig:initial_data}
\end{subfigure}
\hfill
\begin{subfigure}[t]{.19\textwidth}
 \centering
 \includegraphics[height=5cm]{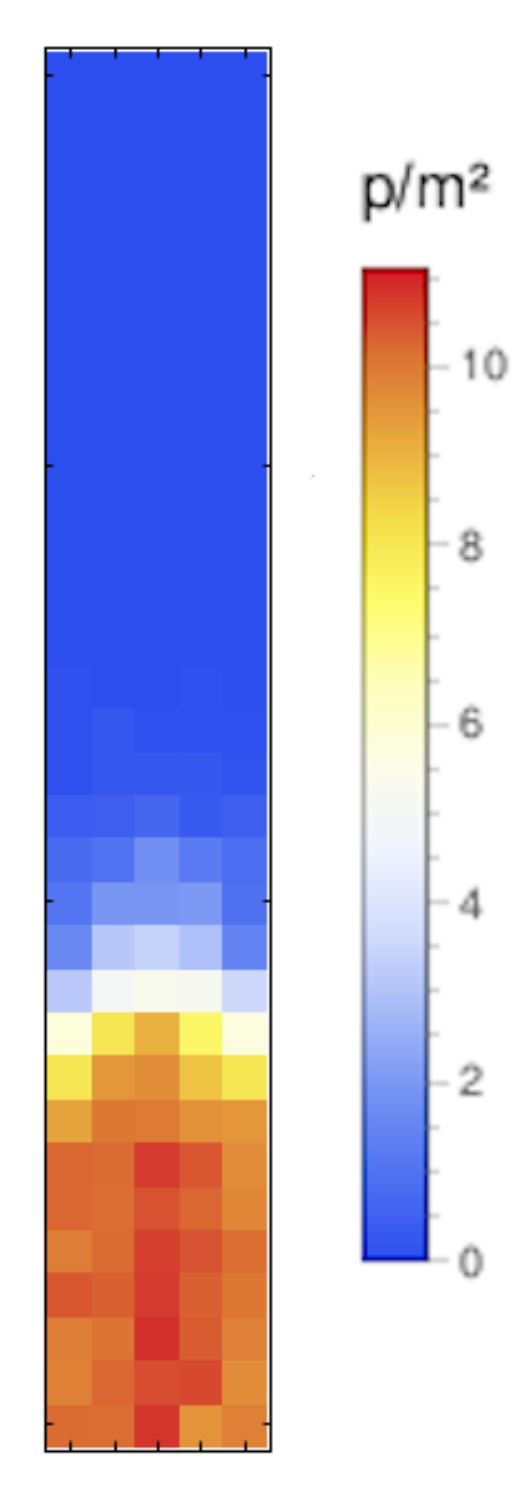}
 \caption{}
 \label{fig:density}
\end{subfigure}
\hfill
\begin{subfigure}[t]{.475\textwidth}
 \centering
 \includegraphics[width=\linewidth]{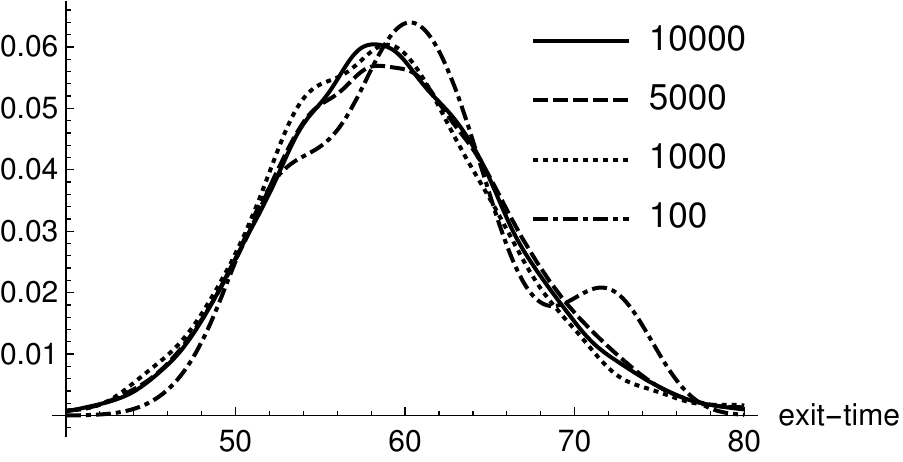}
 \caption{}
 \label{fig:binom}
\end{subfigure}
\caption{
    (A) distribution of $n=40$ agents at time $t=0$ and $t=16 \Delta t$,
    (B) the simulated density,
    (C) densities of exit times when increasing the number of Monte-Carlo runs.
}
\end{figure}

 We check the consistency of the estimated average exit time by varying the number of Monte-Carlo simulations. We observe that the distribution of the exit time
converges to a unimodal curve, see Figure~\ref{fig:binom}. Similar results are obtained across a large range of parameter combinations.

\subsection{Calibration}
In this section, we discuss a possible calibration of the developed CA approach using the experimental data available, see~\cite{hp:juelich}.
We wish to identify the parameter scaling parameter $\beta$, the timestep $\Delta t$ and the exit rate $p_{ex}$. To do so, we make the following assumptions:
\begin{itemize}[leftmargin=10pt,topsep=5pt]
 \item The outflow rate $p_{ex}$ does not depend on the motivation level and the corridor width.
 \item There is a one-to-one relation between the parameter $\beta$ and the time step $\Delta t$.
 \end{itemize}

\noindent We start by considering the dynamics of a single agent in the corridor. These dynamics, although not including any interactions, give first insights and provide reference values for the calibration.

 Velocities of pedestrians are often assumed to be Gaussian distributed. Different values for the mean and variance can be found in the literature, see for example ~\cite{a:speed1, a:speed5}.
We set the desired maximum velocity of a single agent to $v_{max}=1.2\frac{\mathrm{m}}{\mathrm{s}}$, as for example in~\cite{a:speed1}.
Hence a single motivated agent, having motivation level $\mu=1$, needs approximately $8$ seconds to travel the $9.6$m long corridor.

 Let $\bar{N}$ denote the average number of time steps to the exit. We will see in the following that there is a one to one relationship between the scaling parameter
$\beta$ and the exit time, which allows us to estimate the time step $\Delta t$.

\begin{figure}[htb]
 \centering
\begin{subfigure}[t]{0.475\textwidth}
\centering
\includegraphics[height=5cm]{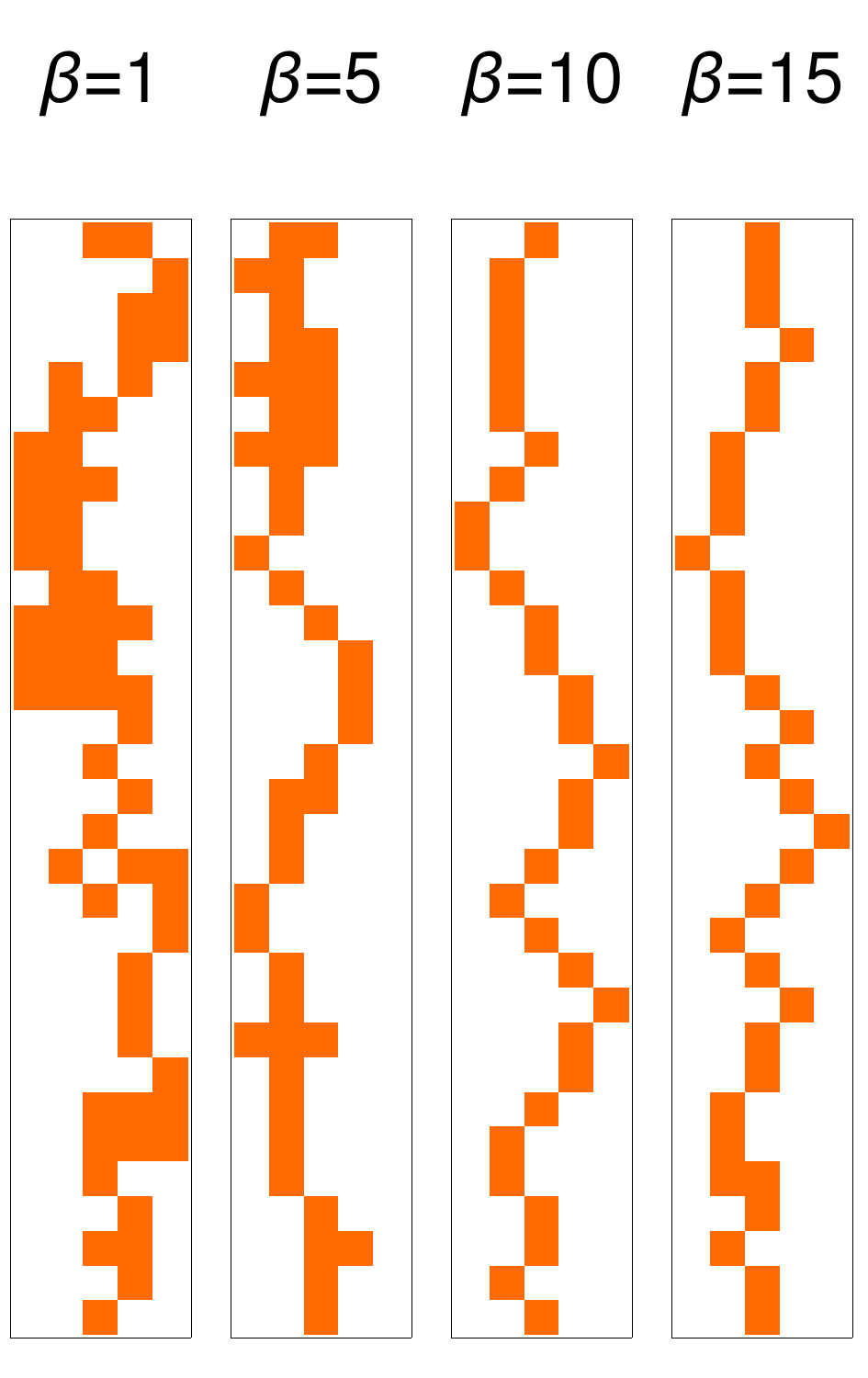}
 \caption{Trajectories for different $\beta$ - increasing $\beta$ reduces the randomness of the walk.}
\label{fig:3betas}
 \end{subfigure}
\hfill
\begin{subfigure}[t]{.475\textwidth}
 \centering
 \includegraphics[width=\linewidth]{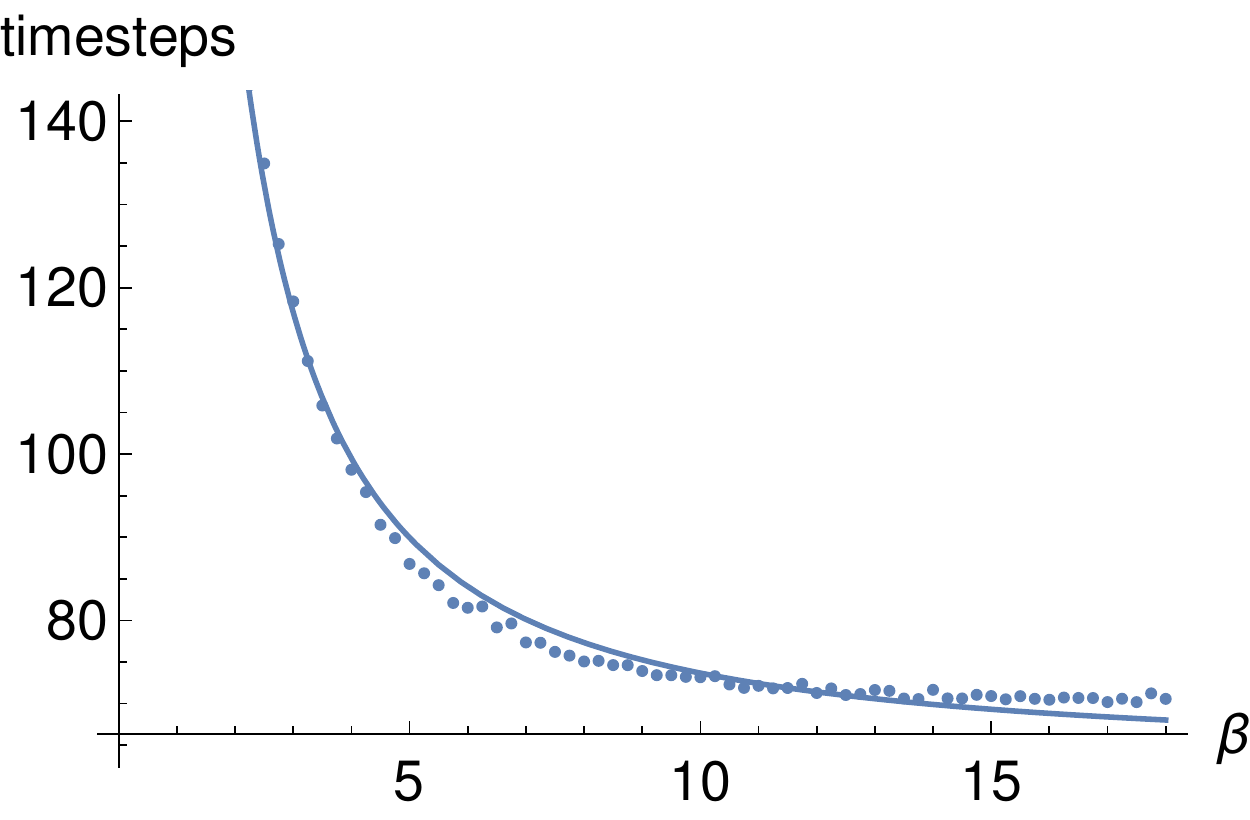}
 \caption{ $\bar{N}$ for different values of $\beta$. The dots mark experimental data, the curve is that of the relation~\eqref{eq:Nbar}.}
 \label{fig:betafunktion}
 \end{subfigure}
\caption{Influence of the scaling parameter $\beta$ on individual dynamics.}
\end{figure}

Figure~\ref{fig:3betas} illustrates the dynamics of this single agent for different values of $\beta$ -- we see that the larger $\beta$, the straighter the path to the exit.
We observe that the average number of time steps $\bar{N}$ to the exit of an agent starting at the same position converges as $\beta$ increases, see Figure~\ref{fig:betafunktion}.
 The observed relation between the exit time and the value of $\beta$ in Figure~\ref{fig:betafunktion} can be estimated by a function of the form
\begin{align}\label{eq:Nbar}
\bar{N}(\beta,p_{ex} =1.1
)=63.528 + \frac{244.082}{\beta^{1.38148}},
\end{align}
which was computed using a least square-approach for $a+\frac{b}{\beta^c}$. The functional relation captures the asymptotic behavior correctly (converging to
the minimum number of steps going straight to the exit) and the sharp increase for small $\beta$.

 This asymptotic relation allows us to estimate the time steps $\Delta t$ for a given value of $\beta$ in case of a single agent.
Since a motivated agent moves on average every second step, it needs approximately $64$ steps to exit the corridor, which corresponds to $32$ vertical fields.
A somehow similar approach was proposed in~\cite{a:cavel}, where the position of agents was updated according to the individual velocity.

 We will now estimate the missing two parameters $\beta$ and $p_{ex}$ using three different data sets, see Table~\ref{table:dataset}.
We restrict ourselves to these three datasets, since the number of individuals in each run is similar and their initial distribution is close to uniform, fitting
the initial conditions of the CA simulations best. For each run, we use the respective modified corridor width $w$, to ensure a consistent discretization of the exit
and the number of agents $n$ as detailed in Table~\ref{table:dataset}.

\begin{table}
    \begin{center}
\begin{tabular}{lcc}
\toprule
Run & $\mu=1$ & $\mu_0$ \\
\midrule
01, $n=1$ & $8\mathrm{s}$ & \\
02, $n=63, w=1.2\mathrm{m}$ & $53\mathrm{s}$ & $64\mathrm{s}$ \\
03, $n=67, w=3.4\mathrm{m}$ & $60\mathrm{s}$ & $68\mathrm{s}$ \\
04, $n=57, w=5.6\mathrm{m}$ & $55\mathrm{s}$ & $57\mathrm{s}$ \\
\bottomrule
\end{tabular}
\end{center}
\caption{exit times for different runs and different motivations from~\cite{art:armin}.
Run 01 is used to set the desired maximum velocity $v_{max}$.}
\label{table:dataset}
\end{table}

\noindent \textit{Reference values:} We use the experimental data to obtain reference values for $\beta$ and $p_{ex}$. For $p_{ex}$, we use
all data sets available, that is a total number of $980$ trajectories recorded for corridors of different widths and consider the respective exit times.
This gives a first approximation $p_{ex} = 1.1\frac{\mathrm{p}}{\mathrm{s}}$, which we use as a reference value for the calibration later on. A similar value for
$p_{ex}$ was reported in~\cite{mt:ben}. We will allow for
estimates within a $50\%$ deviation from that value. Furthermore, we restrict
$\beta$ to $[0.5,10]$ (motivated by the observations in Figure~\ref{fig:3betas}).

 The calibration is then based on minimizing the difference between the observed exit time and the computed average exit time $\bar{T}$. We define the average exit time
\begin{equation*}
\bar T = \bar{T}(\beta, p_{ex}, \Delta t, \mu, n, w): [0,\infty )\times \mathbb{R}^+\times \mathbb{R}^+\times (-\infty,1]\times \mathbb{N}\times \{0.9, 3.3, 5.7\}\rightarrow\mathbb{R}^+\,,
\end{equation*}
that is the time needed for the last agent to leave the domain, for $n$ individuals in a corridor of width $w$ and parameters
 $\beta$, $p_{ex}$, $\Delta t$ and $\mu$.
The calibration is then based on minimizing the functional
\begin{align}\label{eq:minimise2}
\begin{aligned}
\mathcal Z = & \Bigl( (\bar T (\beta, p_{ex}, 63, \Delta t , 0.9)-53)^2+
(\bar T (\beta, p_{ex}, 67, \Delta t , 3.3)-60)^2\\
& + (\bar T (\beta, p_{ex}, 57, \Delta t , 5.7)-55)^2 \Bigr) ^{0.5},
\end{aligned}
\end{align}
using the data stated in Table~\ref{table:dataset}.

 The functional $\mathcal{Z}$ is not differentiable, hence we used derivative free methods to find a minimum.
We first used a parallel Nelder-Mead, which did not converge. We believe that this is caused by the stochasticity of the problem (since we average over $5000$
Monte-Carlo runs to compute the average exit time) as well as the form of the functional itself. Similar problems were reported in ~\cite{a:cal3}.
Systematic computational experiments show that the parameter $\beta$ has a small influence on the exit time. In narrow corridors, increasing the value of $\beta$ does not improve the exit time, since the geometry restricts the range of jumps. In wider corridors, $\beta$ plays a more important role. However, we
have seen that the exit time for a single agent converges as $\beta$ increases. Therefore, we can not expect a unique single optimal value. Furthermore, we believe that the parameter $\beta$ has a smaller influence the more agents are in the corridor.

\begin{figure}[htb]
\scalebox{0.96}{ \centering
 \begin{subfigure}[t]{.475\textwidth}
 \includegraphics[height=4cm]{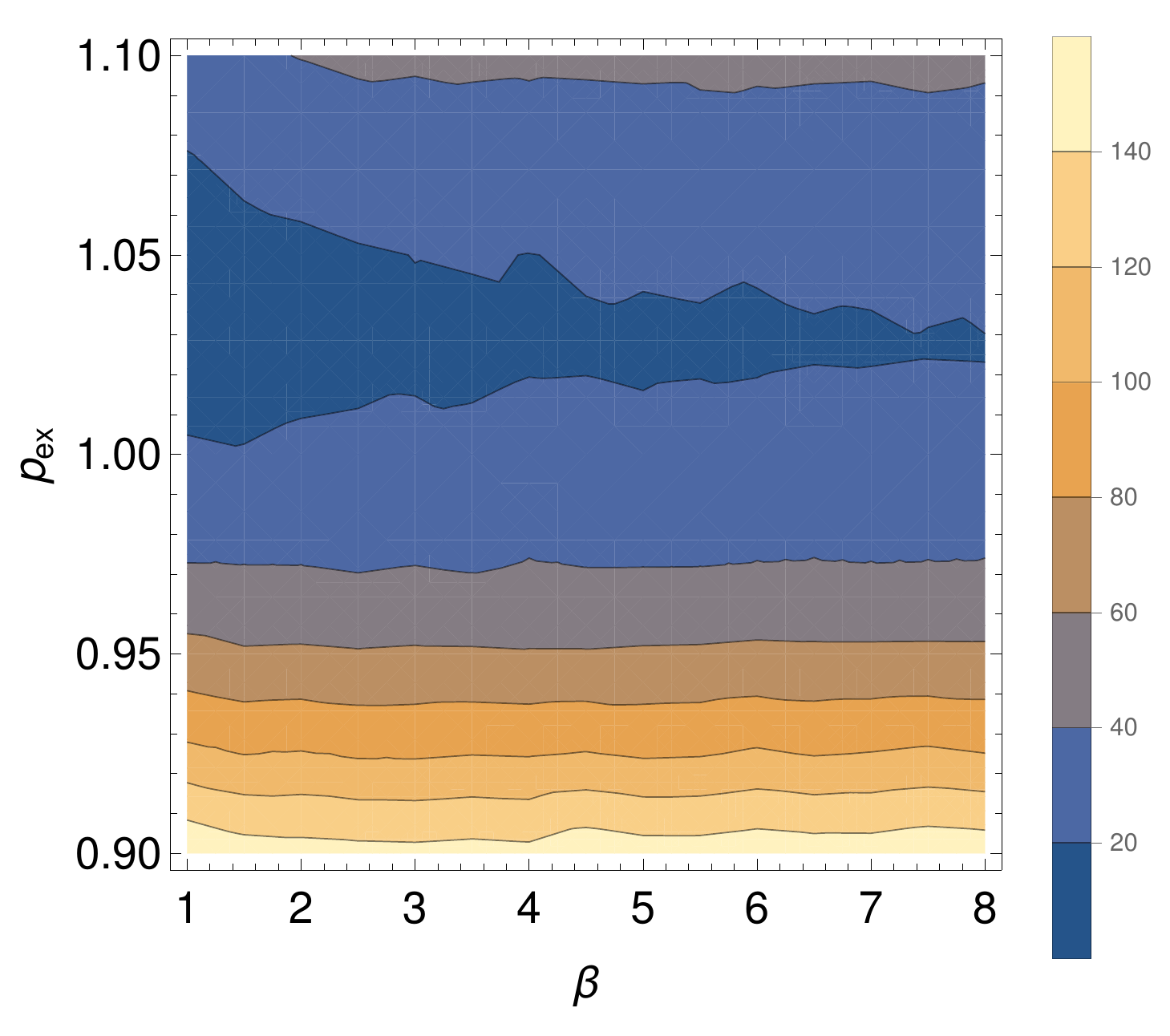}
\caption{The average exit time for a discrete set of $\beta / p_{ex}$-combinations.}
 \label{fig:pexbeta}
\end{subfigure}
\hfill
 \begin{subfigure}[t]{.475\textwidth}
 \includegraphics[height=4cm]{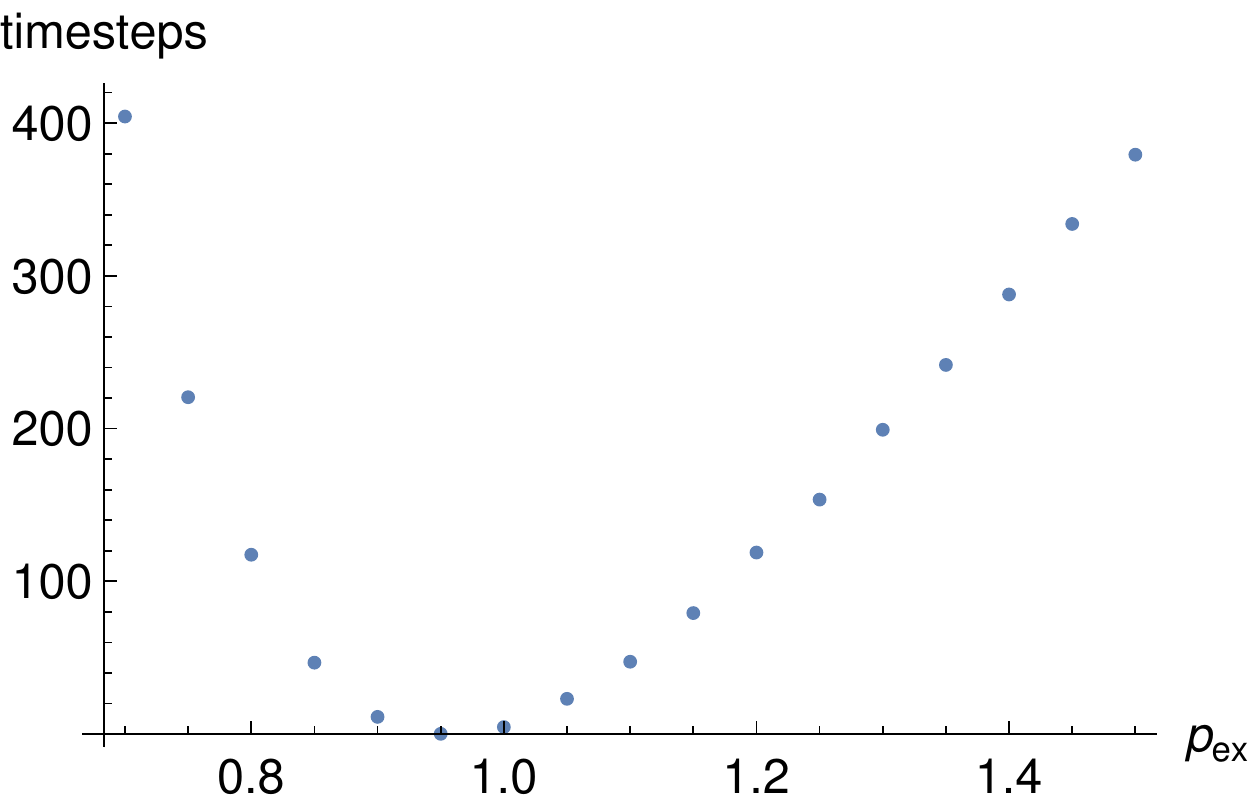}
 \caption{$\bar{T}$ as a function of $p_{ex}$ for a fixed value of $\beta$.}
 \label{fig:guessingpex}
\end{subfigure}}
\caption{Average exit time as a function of $\beta$ and $p_{ex}$, or $p_{ex}$ only.
}
\end{figure}

Finally, we estimate the two parameters $\beta$ and $p_{ex}$ using a discrete search in the range $[0.5,10]\times [0.55,1.65]$. In doing so, we see that the outflow parameter $p_{ex}$ can be clearly estimated for a fixed value of $\beta$, see Figure~\ref{fig:guessingpex}. However, the parameter $\beta$ is much more difficult to determine, as Figure~\ref{fig:pexbeta} shows.
Using the three data sets stated in Table~\ref{table:dataset} we obtain the best fit using
\begin{align}\label{eq:betapex}
(\beta^{\operatorname{min}}, p_{ex}^{\operatorname{min}})\simeq(3.84, 1.15)\,,
\end{align}
which leads to a deviation of $1.04$ seconds in Equation~\eqref{eq:minimise2}.
With a similar approach we then estimate the parameter $\mu_0\simeq-1.22$ for less motivated agents according to Table~\ref{table:dataset}. This results in a speed of $0.57\frac{\mathrm{m}}{\mathrm{s}}$.

\begin{remark}
At first glance, the value $\beta=3.84$ may seem too small given the simulation results shown in Figure \ref{fig:3betas}.
However, Maury suggests a similar value in \cite{book:maury} - in particular $\beta\approx 1/(\Delta x)$ where $\Delta x$ is the cell size.
In our setting this would correspond to the value $3.33$, which is close to the value obtained through from the calibration.
This can be explained by the fact that the effect of $\beta$ is smaller in crowded rooms.
\end{remark}

\subsection{Microscopic simulations}
\begin{figure}[htb]
\scalebox{0.96}{ \centering
\begin{subfigure}[t]{0.44\textwidth}
    \includegraphics[height=4cm]{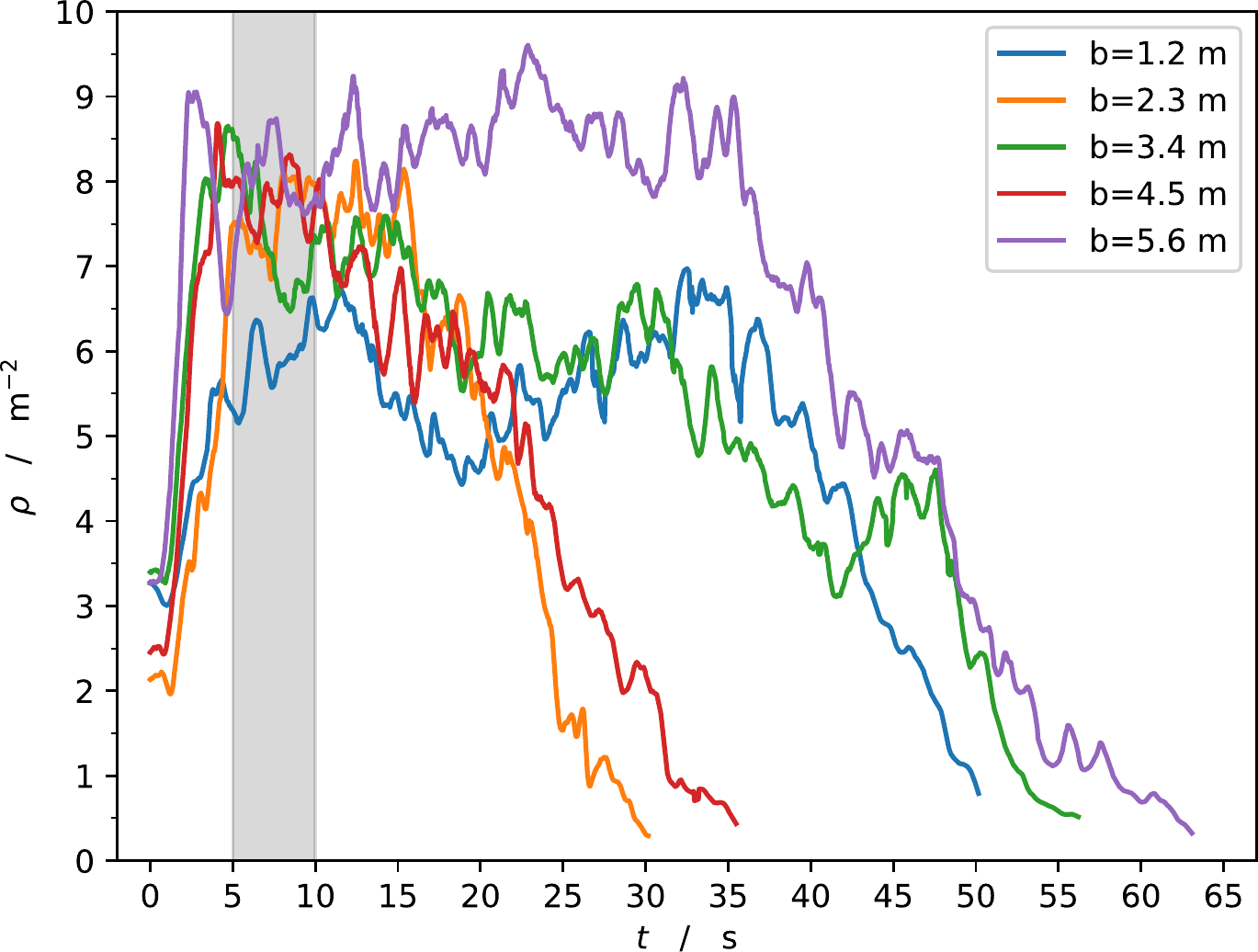}
 \caption{Experimental results, reproduced with permission from~\cite{art:armin}.}
\label{fig:juelichhoch}
 \end{subfigure}
\hfill
\begin{subfigure}[t]{.45\textwidth}
 \includegraphics[height=4cm]{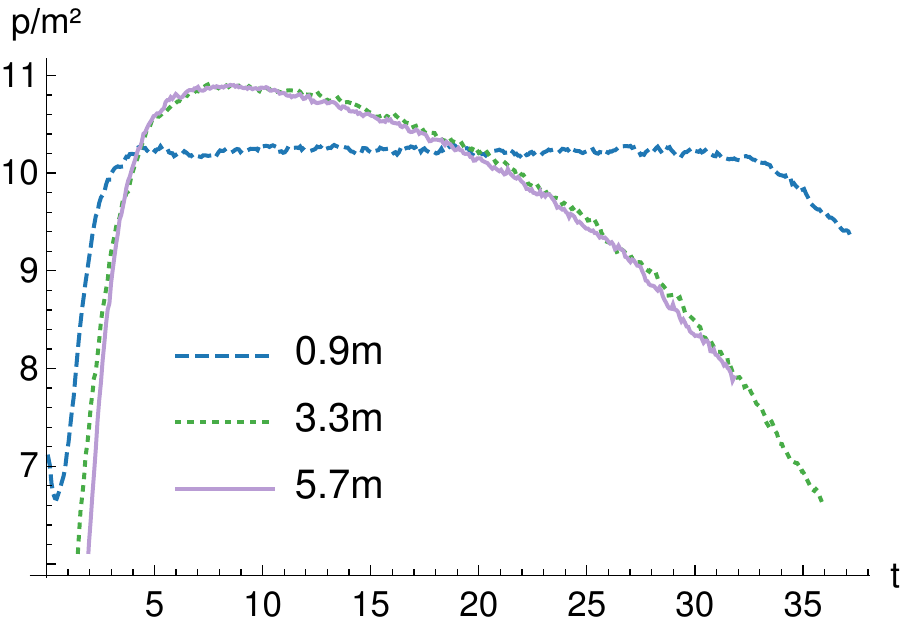}
 \caption{Microscopic simulations for the densities for $n=60$ and $\mu=1$.}
 \label{fig:mikro1}
 \end{subfigure}}
 \caption{Impact of the corridor width on the maximum density. The CA approach yields comparable results for high density regimes and low motivation level.}
\end{figure}

\begin{figure}[htb]
 \scalebox{0.96}{\centering
\begin{subfigure}[t]{0.44\textwidth}
\includegraphics[height=4cm]{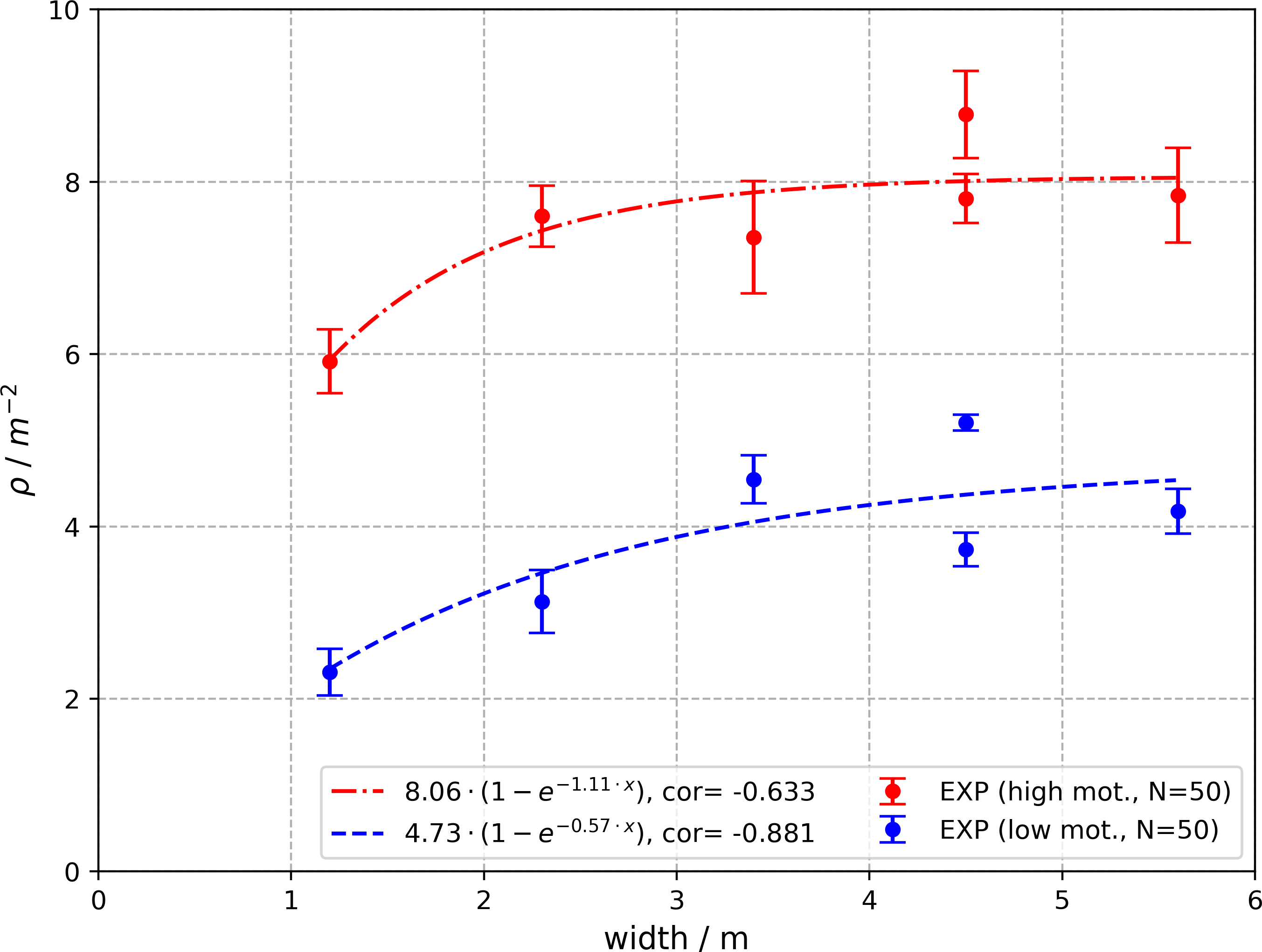}
 \caption{Experimental results (with permission from~\cite{mt:ben})}
\label{fig:ben1}
 \end{subfigure}
\hfill
\begin{subfigure}[t]{.5\textwidth}
    \includegraphics[height=4cm]{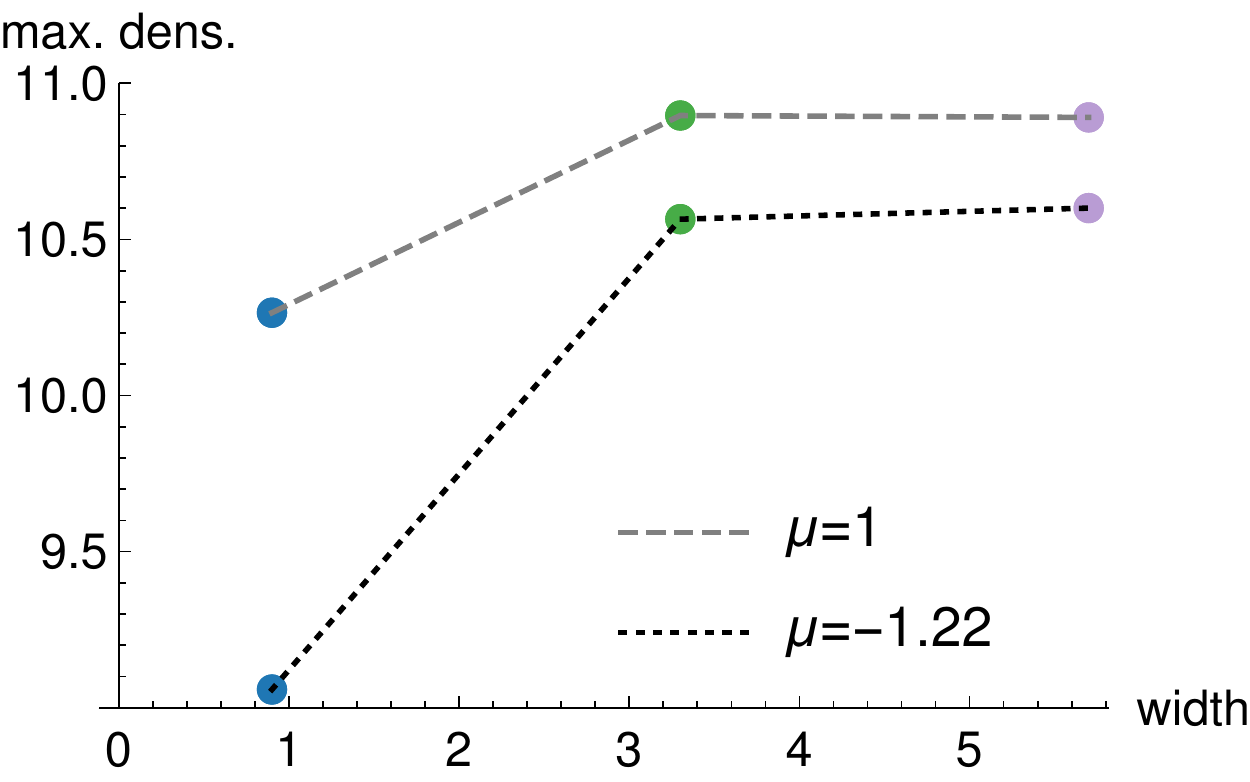}
 \caption{Maximum agent density using the CA model with $n=60$.}
 \label{fig:mudens}
 \end{subfigure}}
 \caption{Impact of the motivation level on the maximum pedestrian density: experimental (A) vs. microscopic simulations (B).}
\end{figure}

We conclude this section by presenting calibrated CA simulations, that are consistent with the experimental data. We observe that wider corridors lead to a higher maximum density in front of the exit
for different motivation levels, see Figure~\ref{fig:mikro1}. Note that this is also the case when changing the outflow rate $p_{ex}$. Higher motivation levels $\mu$ lead
to higher densities in front of the exit as can be seen in Figure~\ref{fig:mudens}. A similar behavior was observed in the experimental results as well as the computational experiments
discussed in~\cite{art:armin, mt:ben}.

\begin{remark}
Note that we observe similar results if we replace the exponential function in~\eqref{eq:transitionrate} by
$\max (0,\phi (x,y)-\phi (x+\Delta x,y))$. However, this function does not satisfy the necessary regularity to
at least formally derive the corresponding macroscopic PDE model.
 \end{remark}

\section{The macroscopic model}\label{s:pde}

\noindent In this section, we derive and study the corresponding macroscopic PDE model, in particular existence of solutions as well
as different options to calculate the path to the exit.

The corresponding macroscopic PDE can be formally derived from the cellular automaton approach discussed in Section~\ref{sec:ca}.
Here, we use a Taylor expansion to develop the transition rates and functions in $x\pm\Delta x$ and $y \pm \Delta x$. This rather tedious calculation
can be done in a systematic manner using a similar approach as discussed in~\cite{KRRW2015}.

\subsection{The PDE and its analysis}
\label{sec:pde analysis}
 We recall that $\rho = \rho(x,y,t)$ denotes the density of pedestrians at position $(x,y)$ and time $t$ and $\phi = \phi(x,y)$ is the potential leading towards the exit $\Gamma_E$.
Let $\Omega \subset \mathbb{R}^2$ denote the domain, $\Gamma_W$ the walls and $\Gamma_E$ the exit with $\Gamma_W \cup \Gamma_E = \partial \Omega$ and $\Gamma_W \cap \Gamma_E = \emptyset$.

Then, the pedestrian density $\rho = \rho(x,y,t)$ satisfies a nonlinear Fokker-Planck equation for all $(x,y) \in \Omega$:
\begin{subequations}\label{e:fpe}
\begin{align}\label{eq:pde}
\partial_t \rho(x,y,t) &= \alpha_\mu \Div \left(\nabla \rho(x,y,t) + 2\beta \rho(x,y,t)(1-\rho(x,y,t))\nabla \phi(x,y)\right)\\
\rho(x,y,0)&=\rho_0(x,y).
\end{align}
The parameter
\begin{align}\label{eq:alpha}
\alpha_\mu:=\frac{1}{8(3-\mu)}
\end{align}
depends on the motivation $\mu$, while $\beta$ corresponds to the ratio between the drift and diffusion.
The function $\rho_0 = \rho_0(x,y)$ is the initial distribution of agents. Equation~\eqref{e:fpe} is supplemented with the following boundary conditions:
\begin{align}\label{eq:bdc}
\begin{aligned}
\textbf{j}\cdot \textbf{n}&=0, \phantom{p_{ex}}\quad \operatorname{ on } \Gamma_W,\\
\textbf{j}\cdot \textbf{n}&=p_{ex}\rho, \quad \, \operatorname{ on } \Gamma_E,
\end{aligned}
\end{align}
\end{subequations}
where $\mathbf{j} = \nabla \rho + 2 \beta \rho (1-\rho) \nabla \phi$ and $\mathbf{n}$ is the unit outer normal vector. We recall that the parameter $p_{ex}$ is
the outflow rate at the exit $\Gamma_E$.

\begin{remark}
Note that the motivation parameter $\mu$ enters the PDE via $\alpha_\mu$ only. It corresponds to a rescaling in time, accelerating
or decelerating the dynamics.
 \end{remark}

 First, we discuss existence and uniqueness of solutions to~\eqref{e:fpe}. Stationary solutions of a similar model were recently investigated by Burger and Pietschmann,
see~\cite{art:burger}. The existence of the respective transient solutions was then shown in~\cite{GSW2019}.
It is guaranteed under the following assumptions:
\begin{enumerate}[label=(A\arabic*)]
\item Let $\Omega\subset \mathbb{R}^2$ with boundary $\partial \Omega $ in $C^2$. \label{a:regdomain}
\item Let $p_{ex}$ be in $[0,1]$. \label{a:outflow}
\item Let $\phi$ be in $H^1(\Omega )$. \label{a:potential1}
\end{enumerate}
Note that assumption~\ref{a:regdomain} is not satisfied in the case of a corridor. However, as pointed out in~\cite{art:burger}, this condition could be relaxed to
Lipschitz boundaries with some technical effort.

\begin{theorem}{(Existence of weak solutions)}\label{theo:weak}
Let assumptions~\ref{a:regdomain}-\ref{a:potential1} be satisfied. Let $\mathcal{S} = \lbrace \rho \in L^2(\Omega): 0 \leq \rho \leq 1\rbrace$ and the initial datum
$\rho_0:\Omega \rightarrow \mathcal S^o$ be a measurable function such that $E(\rho_0) < \infty$, where entropy $E$ is defined by
\begin{align}\label{eq:entropyfunctional}
E(\rho) = \int_{\Omega} \left[\rho \log \rho +(1-\rho) \log (1-\rho)+2\beta \rho \phi\right]dx\,.
\end{align}
Then there exists a weak solution to system~\eqref{e:fpe} in the sense of
\begin{multline}\label{eq:weakpde}
\int_0^T \Bigl[ \langle \partial_t \rho ,\varphi \rangle_{H^{-1},H^1}ds ~- \\
\alpha_\mu\int_\Omega ( (2\beta \rho(1-\rho)\nabla \phi +\nabla \rho))\nabla\varphi dx+p_{ex}\int_{\Gamma_E}\rho \varphi ds \Bigr] dt=0,
\end{multline}
for test functions $\varphi \in H^1(\Omega )$. Furthermore
\begin{align*}
& \partial_t \rho \in L^2 (0,T;H(\Omega)^{-1}),\\
& \rho \in L^2 (0,T;H^1(\Omega)).
\end{align*}
\end{theorem}

 The existence proof is based on the formulation of the equation in entropy variables, that is
\begin{align}\label{eq:entropieform}
\partial_t \rho(x,t) = \Div (m(\rho)\nabla u(x,t)),
\end{align}
where $m(\rho ) = \rho(1-\rho)$ is the mobility function and $u = \frac{\delta E}{\delta \rho} = (\log \rho - \log (1-\rho)+2 \beta \phi)$ the so-called entropy variable.
Note that the proof is a straightforward adaptation of the one presented in~\cite{GSW2019}, hence we omit its details in the following.

\subsection{Moving towards the exit}

In the following we discuss different possible choices for the potential $\phi$.

\subsubsection*{The eikonal equation}
The shortest path to a target, such as the exit $\Gamma_E$, can be computed by solving the eikonal equation, see \cite{book:tosin}:
\begin{equation}\label{eq:eik}
\begin{aligned}
    \| \nabla \phi_E(x,y) \|^2 &=1, &&\textrm{for }(x,y) \textrm{ in } \Omega, \\
\phi_E(x,y)&=0, &&\textrm{on } (x,y)\textrm{ in }\Gamma_E.
\end{aligned}
\end{equation}
Solutions to~\eqref{eq:eik} are in general bounded and continuous, but not differentiable,
see~\cite{BCD2008}. However, in case of the considered corridor geometry, we have the following improved regularity result.
\begin{theorem}\label{theo:reg_eikonal}{(Regularity of $\phi_E$)}
Let $\Omega\subset \mathbb{R}^2$ be a rectangular domain and $\Gamma_E\subset \partial \Omega$ be a line segment in one of the four edges. Then there exists a
solution $\phi_E\in H^1(\Omega)$ to~\eqref{eq:eik}.
\end{theorem}

 The proof can be found in the Appendix and is based on~\cite{BCD2008}, Proposition 2.13.

\subsubsection*{The Laplace equation}\label{chap:laplace}

Alternatively, we consider an idea proposed by Piccoli and Tosin in~\cite{a:Piccoli2011}. Let $\phi_L = \phi(x,y)$ denote the solution of the Laplace equation
on $\Omega \subset \mathbb{R}^2$:
 \begin{equation}\label{eq:laplace1}
 \begin{aligned}
     \Delta \phi_L(x,y)&=0,\phantom{(x,y)} &&\text{for } (x,y)~ \textrm{ in }\Omega, \\
     \phi_L(x,y) &=d(x,y), &&\text{for } (x,y)~ \textrm{ on }\partial \Omega,
\end{aligned}
\end{equation}
where $d = d(x,y)$ corresponds to the Euclidean distance of the boundary points to the exit $\Gamma_E$. Note that in this case of the corridor
the function $d$ is not differentiable at the corners but Lipschitz continuous. Hence standard methods for elliptic equations yield the following
regularity result.

\begin{theorem}\label{theo:reg_laplace}{(Regularity of $\phi_L$)}
Let $d\in C(\partial\Omega )$ defined as above, $\Omega\subset \mathbb{R}^2$ be bounded.
Then there exists a unique solution $\phi_L \in H^1(\Omega)$ to~\eqref{eq:laplace1}.
\end{theorem}

 The proof can be found in ~\cite{book:1968linear}, Section 5.

\medskip
 In the following, we discuss the similarities and difference of the potentials $\phi_E$ and $\phi_L$.
In the case of a 1D corridor with a single exit, that is a line with a single exit on one of the two endpoints, the potentials are identical.
However, in the case of two exits at the respective endpoints, the Laplace equation gives $\phi_L \equiv 0$, which does not provide a sensible potential.

\begin{figure}[htb]
 \centering
\begin{subfigure}[t]{0.49\textwidth}
\centering
 \includegraphics[height=4.5cm]{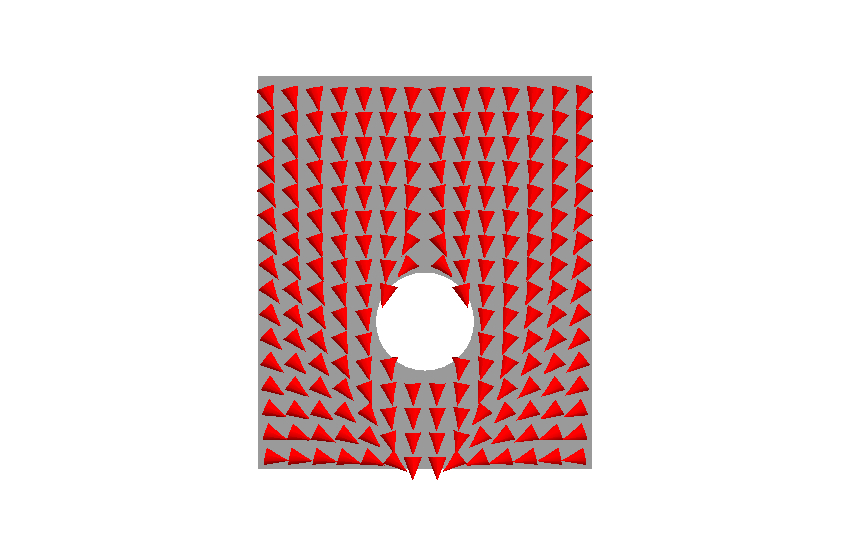}
 \caption{Eikonal equation.}
 \label{fig:eik1}
\end{subfigure}
\hfill	
\begin{subfigure}[t]{0.49\textwidth}
\centering
\includegraphics[height=4.5cm]{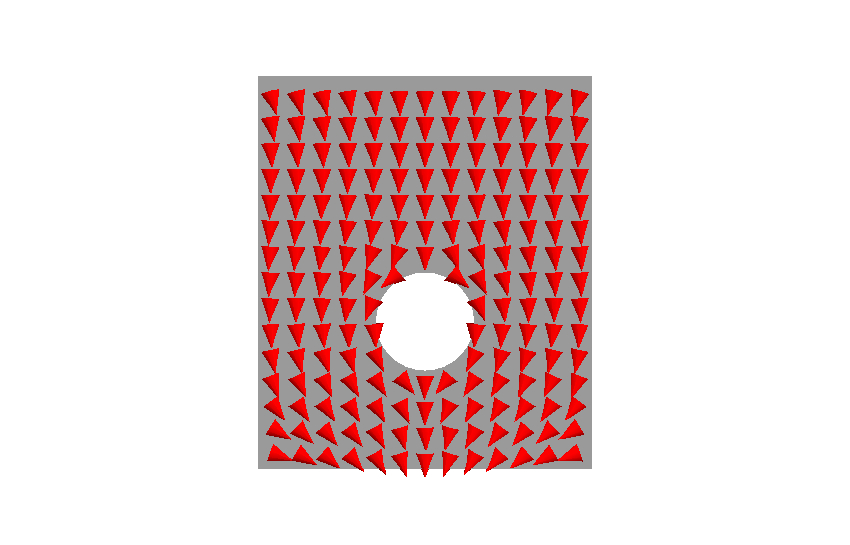}
\caption{Laplace equation (with Neumann bc at the obstacle).}
\label{fig:lap1}
\end{subfigure}
\caption{Comparison of the potentials $\phi_E$ and $\phi_L$ for a convex obstacle.}\label{f:comp_pot_conv}
\end{figure}

\begin{figure}
\begin{subfigure}[t]{0.49\textwidth}
\centering
 \includegraphics[height=4.5cm]{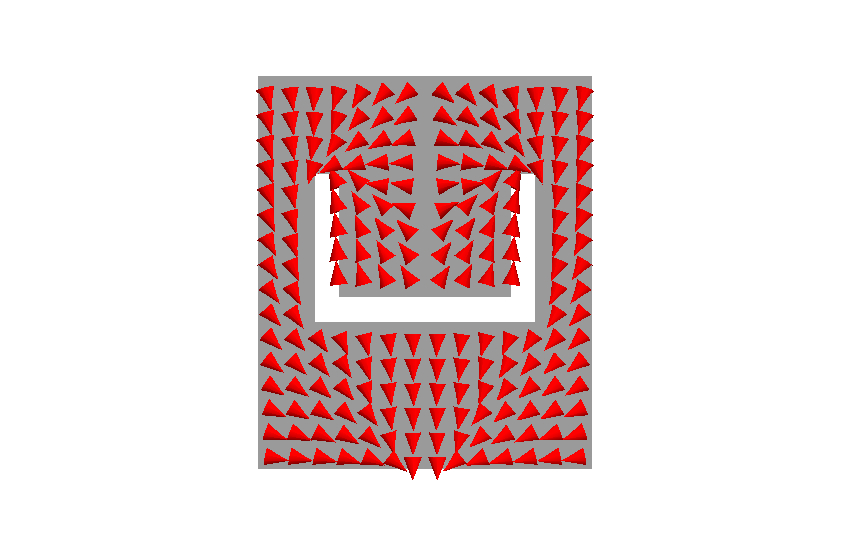}
 \caption{Eikonal equation.}
 \label{fig:eik2}
\end{subfigure}
\hfill
\begin{subfigure}[t]{0.49\textwidth}
\centering
\includegraphics[height=4.5cm]{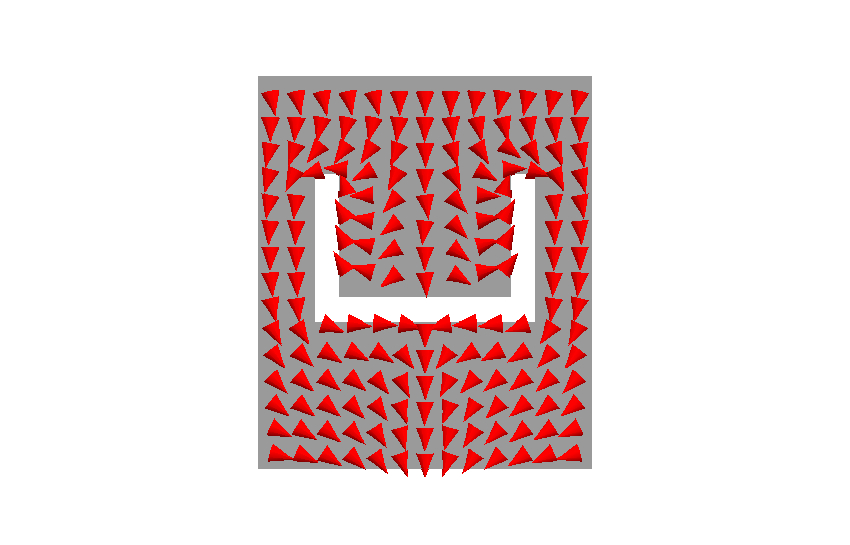}
\caption{Laplace equation (with Neumann bc at the obstacle).}
\label{fig:lap2}
\end{subfigure}
\caption{Comparison of the potentials $\phi_E$ and $\phi_L$ for a non-convex obstacle.}\label{f:comp_pot}
\end{figure}

Figures~\ref{f:comp_pot_conv} and~\ref{f:comp_pot} illustrate the differences between $\phi_L$ and $\phi_E$ in 2D. Note that we choose homogeneous Neumann boundary conditions
at the obstacle walls when solving the Laplace equation~\eqref{eq:laplace1}. We observe good agreement in the case of convex obstacles, see Figure~\ref{f:comp_pot_conv}.
In the case of non-convex obstacles, such as the U-shaped obstacle in Figure~\ref{f:comp_pot}, individuals would first get trapped inside the U using the Laplace equation.
Solving the eikonal equation~\eqref{eq:eik} is in general computationally more expensive than the Laplace equation~\eqref{eq:laplace1}. However,
these costs are negligible since the potential is stationary and computed only once.

\subsection{Characteristic calculus}\label{sec:cc}

We now consider the corresponding inviscid macroscopic model, which can be derived using a different scaling limit from the CA approach.
We focus on the one dimensional case only as we can derive solutions explicitly. A similar problem (with different boundary conditions)
was partially analyzed in~\cite{a:hughes1}.

 The inviscid PDE reduces to a scalar conservation law, posed on $\mathbb{R}_+$, of the form
\begin{equation}
 \partial_t \rho + \partial_x j(\rho) = 0\,,
 \label{eq:scalar conservation law}
\end{equation}
where the flux function is $j(\rho) = -\rho(1-\rho)$. Note that this flux corresponds to the potential $\phi(x) = x$, hence individuals move to the left.
We consider the initial condition
\begin{equation*}
 \rho(x,0) = \rho_0\, \chi_{[0, L]}\,,
\end{equation*}
for some positive $L$, where $\chi$ denotes the characteristic function.
At the origin, we wish to enforce a similar outflow condition as in the viscous case and set $j(0,t) = p_{ex}\, \rho,~t>0$. This
is equivalent to the Dirichlet boundary condition $\rho(0, t) = 1 - p_{ex}$ for all times $t>0$, where we recall that $0<p_{ex} \le 1$. This is an ill-posed problem in general~\cite{leroux_approximation_1979}, and the boundary condition must be relaxed.

 Away from discontinuities, the speed of characteristics is given by
\begin{equation*}
 j'(\rho) = -(1-2\rho)\,.
\end{equation*}
We see that they either point in- or outside of the domain, depending on the magnitude of $\rho$. Recall that for a shock located at $s(t)$, the Rankine-Hugoniot condition reads $\llbracket j'(\rho) \rrbracket = \dot{s}(t) \llbracket \rho \rrbracket$, where $\llbracket f \rrbracket = f^- - f^+$, with $f^\pm(x) = f(x\pm 0)$.
For our choice for $\rho_0$, there is an initial shock at $x_r=L$, which is moving (left) at a speed of
\begin{equation*}
 \dot{x}_r
 = -(1-\rho_0)\,.
\end{equation*}
The larger the initial pedestrian density, the slower the shock moves or the people get closer to the exit.
One can easily check that such a profile satisfies the so-called \emph{Lax entropy condition}, since
\begin{equation*}
 -1 = j'(0) \le \dot{x}_r \le j'(\rho_0)\,,
\end{equation*}
is it therefore admissible.

\begin{figure}
 \begin{center}
 \begin{tikzpicture}[scale=5.1]
 \draw (0,0) -- (1,0) -- (1,1) -- (0,1) -- (0,0);
 \draw[very thick] (0,0) -- (0.5,0.5) -- (1,0.5);
 \draw[very thick] (0,1) -- (0.5, 0.5);

 \node (0) at (-0.05,-0.05) {$0$};
 \node (0one) at (-0.05,1) {$1$};
 \node (one0) at (1,-0.05) {$1$};

 \node (labelx) at (0.5,-0.05) {$p_{ex}$};
 \node (labely) at (-0.05,0.5) {$\rho_0$};

 % solution remains constant
 \coordinate (A) at (0.6,0.15);
 \draw ($ (A) - (0.05,0) $) -- ($ (A) + (0.3,0) $ );
 \draw ($ (A) - (0,0.05) $) -- ($ (A) + (0,0.2) $ );
 \draw[thick] ($ (A) + (0,0.15) $) -- ($ (A) + (0.25,0.15) $ );
 \draw[dotted] ($ (A) + (0.25,0) $) -- ($ (A) + (0.25,0.15) $ );
 \draw[->] ( $(A) + (0.22, 0.075)$ ) -- ( $(A) + (0.15, 0.075)$ );
 \draw ( $(A) + (-0.02,0.15)$ ) -- ( $(A) + (0,0.15)$ );
 \node (rhoA) at ( $(A) + (-0.07,0.15)$ ) {$\rho_0$};

 % right moving front from origin
 \coordinate (B) at (0.07,0.4);
 \draw ($ (B) - (0.05,0) $) -- ($ (B) + (0.3,0) $ );
 \draw ($ (B) - (0,0.05) $) -- ($ (B) + (0,0.2) $ );
 \draw[thick] ($ (B) + (0,0.15) $) -- ($ (B) + (0.10,0.15) $ );
 \draw[dotted] ($ (B) + (0.10,0.08) $) -- ($ (B) + (0.10,0.15) $ );
 \draw[->] ( $(B) + (0.02, 0.12)$ ) -- ( $(B) + (0.08, 0.12)$ );
 \draw[thick] ($ (B) + (0.10,0.08) $) -- ($ (B) + (0.25,0.08) $ );
 \draw[dotted] ($ (B) + (0.25,0) $) -- ($ (B) + (0.25,0.08) $ );
 \draw[->] ( $(B) + (0.22, 0.04)$ ) -- ( $(B) + (0.15, 0.04)$ );
 \node (rhoB) at ( $(B) + (0.3,0.075)$ ) {$\rho_0$};
 \node (pexB) at ( $(B) + (0.23,0.15)$ ) {$1-p_{ex}$};

 % rarefaction wave from origin
 \coordinate (C) at (0.65,0.67);
 \draw ($ (C) - (0.05,0) $) -- ($ (C) + (0.3,0) $ );
 \draw ($ (C) - (0,0.05) $) -- ($ (C) + (0,0.2) $ );
 \draw[thick] ($ (C) + (0,0.075) $) -- ($ (C) + (0.10,0.15) $ );
 \draw[thick] ($ (C) + (0.10,0.15) $) -- ($ (C) + (0.25,0.15) $ );
 \draw[dotted] ($ (C) + (0.25,0) $) -- ($ (C) + (0.25,0.15) $ );
 \draw[->] ( $(C) + (0.22, 0.075)$ ) -- ( $(C) + (0.15, 0.075)$ );
 \draw[dotted] ( $(C) + (0.1,0.15) $ ) -- ( $(C) + (0.1,0.2)$ );
 \draw[->] ( $(C) + (0.1, 0.2)$ ) -- ( $(C) + (0.2, 0.2)$ );
 \draw ( $(C) + (-0.02,0.075)$ ) -- ( $(C) + (0,0.075)$ );
 \node (rhoC) at ( $(C) + (0.3,0.15)$ ) {$\rho_0$};
 \node (maxC) at ( $(C) + (-0.17,0.075)$ ) {$\frac{1}{2}\vee1{-}p_{ex}$};
 \end{tikzpicture}
 \includegraphics[width=0.5\textwidth]{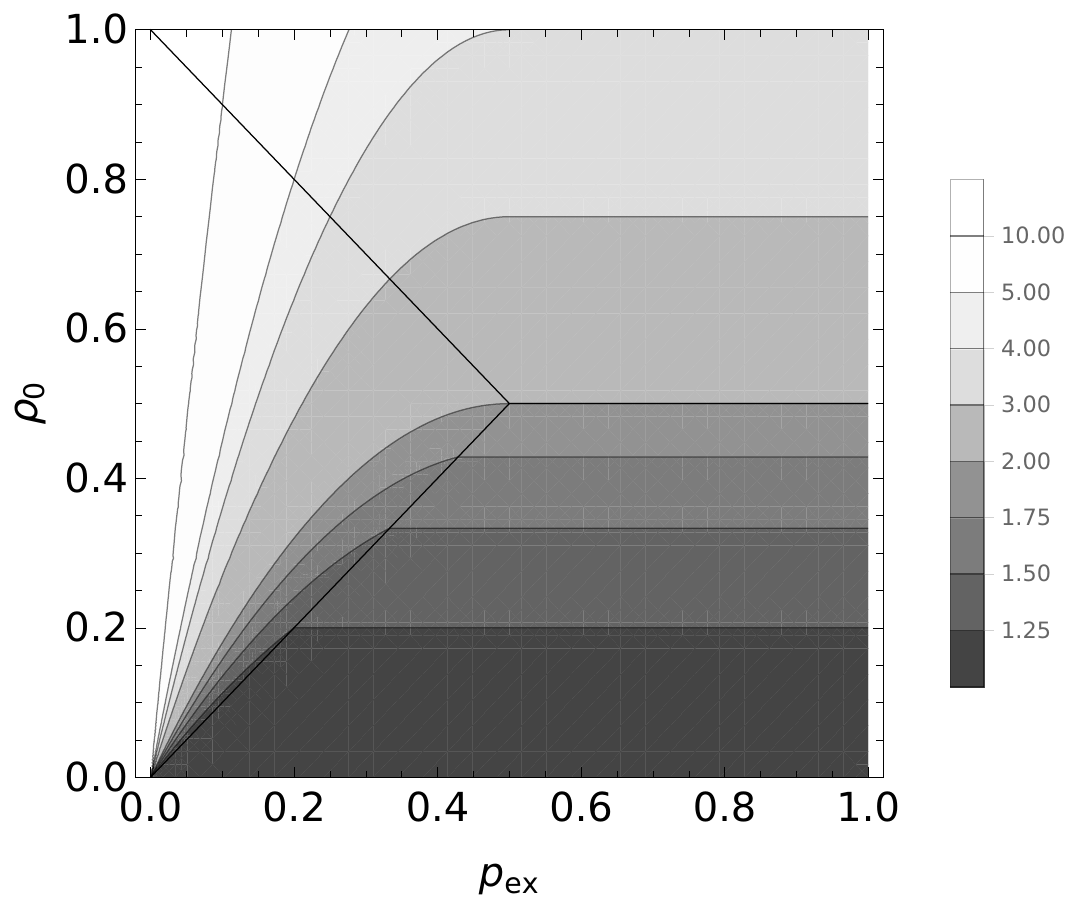}
 \end{center}
 \caption{Left: Bifurcation diagram detailing the behavior of the solution to~\eqref{eq:scalar conservation law}-\eqref{eq:relaxed boundary condition}. The behavior along the interface lines is identical as in the bottom right corner. Right: exit time corresponding to $L=1$.}
 \label{fig:characteristics bifurcation}
\end{figure}
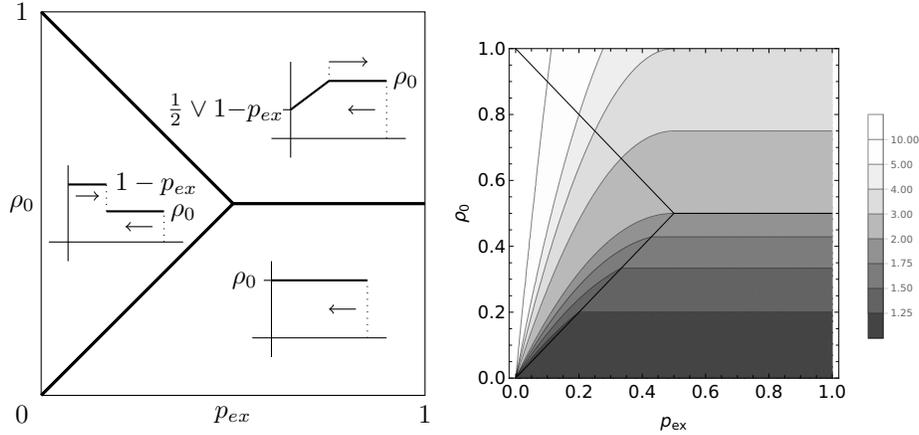

 Next we discuss the behavior of solutions at the exit $x=0$. The proper way to enforce the Dirichlet boundary condition
is derived in~\cite{lefloch_explicit_1988-1, lefloch_explicit_1988}, and reads as follows:
\begin{equation}
 \rho^+(0) \in \mathcal{E}[1-p_{ex}] :=
\begin{cases}
 \left[0, p_{ex}\right] \cup \left\{1-p_{ex}\right\} & \text{ if } p_{ex} < \frac{1}{2}\,,
 \\
 \left[0, \frac{1}{2}\right] & \text{ if } p_{ex} \ge \frac{1}{2}\,.
\end{cases}
\label{eq:relaxed boundary condition}
\end{equation}
Depending on the slope of the characteristics as well as the value of the outflow rate $p_{ex}$, we observe three different cases, which are detailed below and illustrated
in Figure~\ref{fig:characteristics bifurcation}.
\begin{itemize}
 \item \emph{A constant profile} for $\rho_0 \le p_{ex} < \frac{1}{2}$ or $\rho_0 \le \frac{1}{2} \le p_{ex}$. In this case, the characteristics have a negative slope
 and $\rho_0$ is an admissible boundary value. The function $\rho$ vanishes when the shock originating at $x=L$ reaches the origin at time $t = \frac{L}{1-\rho_0}$.
 The situation is similar in the case $\frac{1}{2} < \rho_0 = 1-p_{ex}$, for which characteristics are going inwards but where $\rho_0 \in \mathcal{E}[1-p_{ex}]$. This case
 is illustrated on the lower right in Figure~\ref{fig:characteristics bifurcation}.
 \item \emph{A shock originating at $x=0$} for $p_{ex} < \frac{1}{2}$ and $p_{ex} < \rho_0 < 1-p_{ex}$. In this case, $\rho_0 \notin \mathcal{E}[1-p_{ex}]$, but
 $1-p_{ex} \in\mathcal{E}[1-p_{ex}]$, which we therefore set as a boundary value. This causes a shock at the origin, which travels to the right with speed
 \begin{equation*}
 \dot{x}_l = \frac{-p_{ex}\left(1-p_{ex}\right) + \rho_0\left(1-\rho_0\right)}{1 - p_{ex} - \rho_0}\,,
 \end{equation*}
 until it collides with the back-shock. The collision time and position, $t = t_1^*$ and $x=x_1^*$ respectively, can be calculated from
 $x_1^* = \dot{x}_l t_1^* = L + \dot{x}_r t_1^*$. We obtain
 \begin{equation*}
 t_1^* = \frac{L}{1-p_{ex}}\,, \quad \text{and} \quad x_1^* = \left(1 - \frac{1 - \rho_0}{1-p_{ex}}\right)L\,.
 \end{equation*}
 The resulting shock will then move to the left again, with speed $-p_{ex}$ and reaches the origin at time $t = \frac{\rho_0 \, L}{p_{ex}\left(1-p_{ex}\right)}$. This situation is shown in the
 center left part of Figure~\ref{fig:characteristics bifurcation}.
 \item \emph{A rarefaction wave originating at $x=0$} for $\rho_0 > \frac{1}{2}$ or $\rho_0 > 1-p_{ex}$. In this case, a rarefaction wave will connect the value
 at the boundary, that is $\bar{\rho} = \frac{1}{2} \vee 1-p_{ex}$ , with the state $\rho_0$.
 The rarefaction wave is of the form $\rho(x, t) = \frac{x+t}{2t}$. More precisely we have for any $x > 0$:
 \begin{equation}
 \rho(x,t) =
 \begin{cases}
 \;\bar{\rho} & \text{ if } \frac{x}{t} \le \left(2\bar{\rho}-1\right)\,,
 \\
 \frac{x+t}{2t} & \text{ if } \left(2\bar{\rho}-1\right) < \frac{x}{t} < \left(2\rho_0-1\right)\,,
 \\
 \;\rho_0 & \text{ if } \frac{x}{t} \ge \left(2\rho_0-1\right)\,.
 \end{cases}
 \label{eq:rarefaction wave}
 \end{equation}
 Note that for $1-p_{ex} > \frac{1}{2}$, the constant value $\bar{\rho} = 1-p_{ex}$ is transported into the domain at speed $2\bar{\rho} - 1$.
 The crest the rarefaction wave travels at speed $2\rho_0-1$ until it hits the back-shock at time $t_2^*=\frac{L}{\rho_0}$, for $x_2^* = \left(\frac{2\rho_0-1}{\rho_0}\right)L$.
 This results at a new shock, which originates at position $x_s$ and with velocity $\dot{x}_s = -(1-\rho(x_s))$. From~\eqref{eq:rarefaction wave}
 we also have $\rho(x_s) = \frac{x_s + t}{2t}$. Solving the resulting equation with initial condition $x_s(t_2^*) = x_2^*$ yields
 \begin{equation*}
 x_s(t) = 2\sqrt{L\,\rho_0}\sqrt{t} - t\,.
 \end{equation*}
 If $1-p_{ex} < \frac{1}{2}$, this new shock reaches $0$ at time $t_3^* = 4\rho_0 L$.
 Otherwise, the back-shock will meet the constant state $1-p_{ex}$ at time $t_4^* = \frac{L \rho_0}{(1-p_{ex})^2}$, for $x_4^* = \frac{L \left(1-2p_{ex}\right)\rho_0}{(1-p_{ex})^2}$,
 resulting in a single constant state $\rho(x, t_4^*) = (1-p_{ex})\chi_{[0,x_4^*]}$. This constant profile
 then moves with speed $-p_{ex}$, and reaches the origin at time $t = \frac{\rho_0 \, L}{p_{ex}\left(1-p_{ex}\right)}$, see upper right corner of Figure~\ref{fig:characteristics bifurcation}.
\end{itemize}

 Figure~\ref{fig:characteristics bifurcation} illustrates how the exit time changes with the initial pedestrian density and the outflow rate. We see that for an outflow rate $p_{ex} \gtrsim \frac{1}{2}$,
the initial density $\rho_0$ has a much stronger influence on the exit time compared to the value of $p_{ex}$. The situation is somehow reversed for small $p_{ex}$.

\subsection{Numerical results}
We conclude by presenting computational results on the macroscopic level. All simulations use the finite element library Netgen/NgSolve.

We consider a rectangular domain with a single exit as shown in Figure~\ref{fig:domaindisc}
and discretize it using a triangular mesh of maximum size $h=0.1$. The potential $\phi$ is calculated in a preliminary step, by either solving the
eikonal equation~\eqref{eq:eik} or the Laplace equation~\eqref{eq:laplace1}. We use a fast sweeping scheme for the eikonal equation,
as it can be generalized to triangular meshes, see~\cite{a:zhao2}. The discretization of the nonlinear Fokker-Planck equation~\eqref{e:fpe} is based on a $4^{\text{th}}$
order Runge-Kutta method in time and a hybrid discontinuous Galerkin method in space, see~\cite{lehrenfeld}.

\begin{figure}[htb]
 \centering
 \begin{subfigure}[t]{0.45\textwidth}
\includegraphics[height=4cm]{figure_15}
 \caption{Microscopic simulations.}
\label{fig:macro2b}
\end{subfigure}
\begin{subfigure}[t]{0.45\textwidth}
\includegraphics[height=4cm]{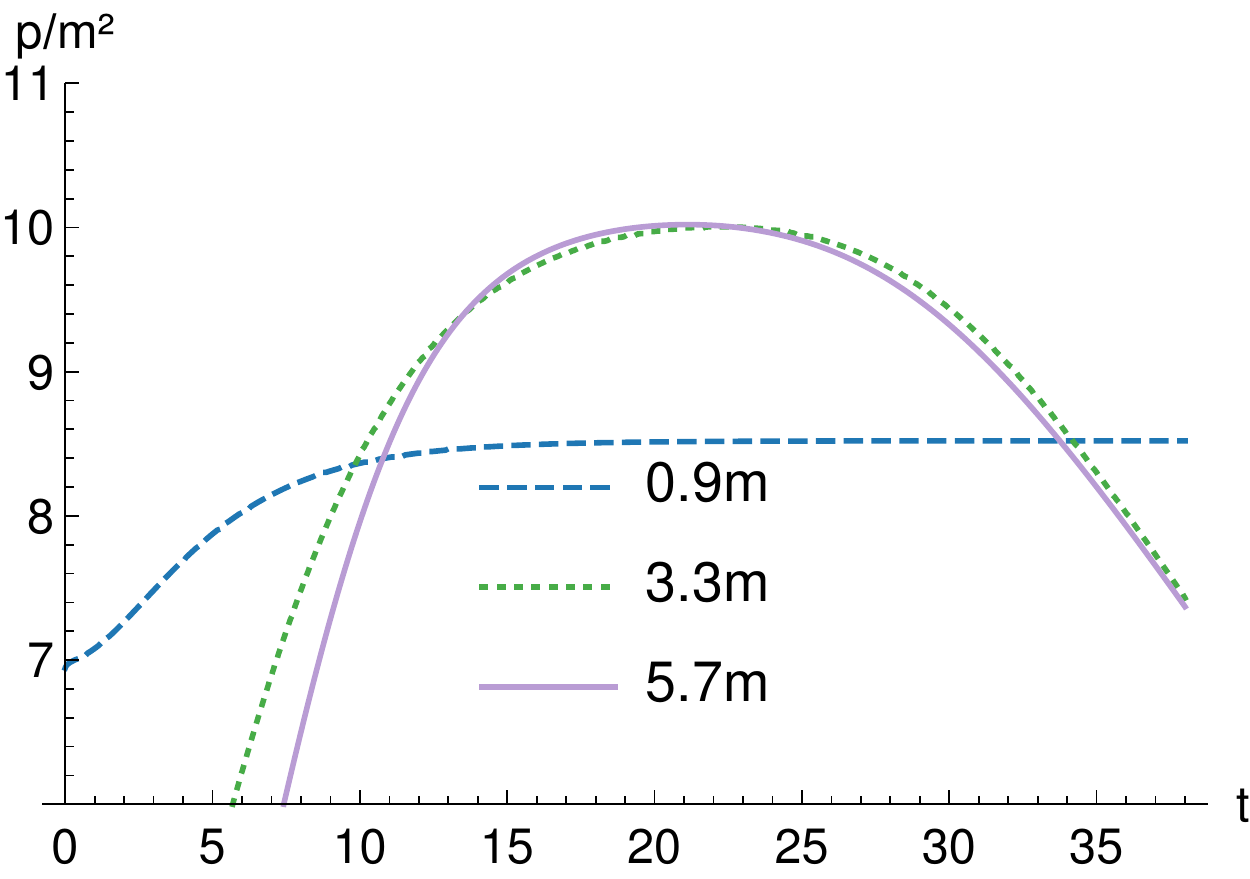}
\caption{Macroscopic simulations.}
\label{fig:macro2a}
\end{subfigure}
\caption{Simulations for $n=60$, $\beta=3.84$, $\mu=1$. We observe a good agreement between the CA (microscopic) and PDE (macroscopic) solutions. The effect of higher densities for wider corridors also occur on a macroscopic scale. There is a clear difference in behavior between narrow and wide corridors}\label{fig:macro1}
\end{figure}

We choose a constant initial datum $\rho_0$, taken such that $\int_\Omega \rho_0 dx=\frac{n}{\rho_s}$. We recall that $\rho_s$ corresponds to the typical pedestrian
density $\rho_s=11.11\frac{\mathrm{p}}{\mathrm{m}^2}$ and $n$ to the number of individuals. The simulation parameters are set to
\begin{align*}
 \beta = 3.84,~ p_{ex} = 1.15,~ \Delta t =10^{-5},\textrm{ and } \alpha_\mu =\frac{1}{16}
\end{align*}
We calculate the densities in the rectangular area highlighted in Figure~\ref{fig:domaindisc}. The macroscopic simulations confirm the microscopic results.
Again, higher densities for wider corridors are observed, see Figure ~\ref{fig:macro1}.

\section{Alternative modeling approaches}\label{s:alt}
We have seen that the proposed CA approach proposed in Section~\ref{sec:ca} reproduces some features of the observed dynamics on the
microscopic as well as on the macroscopic level. In the following, we discuss possible alternatives and
generalizations, which we expect to result in even more realistic results.

\subsection{Density dependent cost}
Hughes~\cite{a:Hughes0} proposed that the cost of moving should be proportional to the local pedestrian density. In particular, moving through
regions of high density is more expensive and therefore less preferential. This corresponds to a density dependent (hence time dependent) right-hand side in~\eqref{eq:eik}.
In particular, Hughes proposed a coupling via
\begin{equation}\label{eikonal_coupled}
\begin{aligned}
\lVert \nabla \phi(x,y,t) \rVert& =\frac{1}{1-\rho(x,y,t)},&& \text{ for } (x,y)\text{ in }\Omega\\
\phi(x,y)& =0,\phantom{\frac{1}{1-\rho(x,y,t)}} && \text{ for } (x,y)\text{ in }\Gamma_E .
\end{aligned}
\end{equation}
\noindent We see that the right hand side, which corresponds to the cost of moving, becomes unbounded as $\rho$ approaches the scaled maximum density $1$.
Such density dependent cost should lead to more realistic dynamics. However the analysis of the coupled problem \eqref{e:fpe}-\eqref{eikonal_coupled}
is open. Solutions to~\eqref{eikonal_coupled} have a much lower regularity than required in Theorem~\ref{theo:reg_eikonal}. We expect this to lead to similar
analytic challenges as reported in~\cite{a:hughes1}.
\subsection{Alternative ways to model motivation}
In the following, we discuss different possibilities to include the influence of the motivation level on the dynamics. First by modifying the transition rates and
second by changing the transition mechanism, allowing for shoving.

\subsubsection*{Alternative transition rates}
In the transition rate \eqref{eq:transitionrate}, the motivation relates to the probability of jumping as detailed in Remark \ref{rem:jumping}. It is therefore directly correlated to the
 agent's velocity on a microscopic level. However, one could assume that the motivation increases the probability to move along the shortest path. This could be modeled by transition rates of the form
\begin{align}
\mathcal T^{ij}(x,y) =
\frac{1}{8}
 \exp (\mu \beta(\phi (x,y)-\phi (x+i\Delta x,y+j\Delta x))).
\end{align}
Then the corresponding macroscopic model reads
\begin{align}\label{pde2}
\partial_t \rho(x,y,t) = \frac{1}{8}\Div \left(\nabla \rho(x,y,t) + 2\mu\beta \rho(x,y,t)(1-\rho(x,y,t))\nabla \phi(x,y)\right).
\end{align}
We see that the motivation level $\mu$ enters only in the convective term. Hence higher motivation is directly correlated to a higher average velocity of the crowd on a macroscopic level.

\subsection{Pushing and shoving}
\subsubsection*{Microscopic modeling}
In the previously proposed model, the transition rates depended on the availability of a site and the motivation level. Another possibility to include the latter
is by allowing individuals to push. Different pushing mechanisms have been proposed in the literature. In local pushing models, individuals are only able
to push one neighbor into an adjacent vacant site, while in global pushing individuals can push a given number of neighbors into a direction.
However, since individuals can induce movements of other individuals in some distance (and not only the neighboring sites), an implementation on bounded
domains is not straight forward. In particular, it is not clear how to adapt boundary conditions in case of global pushing, as considered in~\cite{a:pushing2}.
In contrast, local pushing mechanisms, can be translated one-to-one on bounded domains, see~\cite{a:pushing1}.

\begin{figure}[htb]
\centering
\begin{subfigure}[t]{0.45\textwidth}
\centering
\includegraphics[height=4cm]{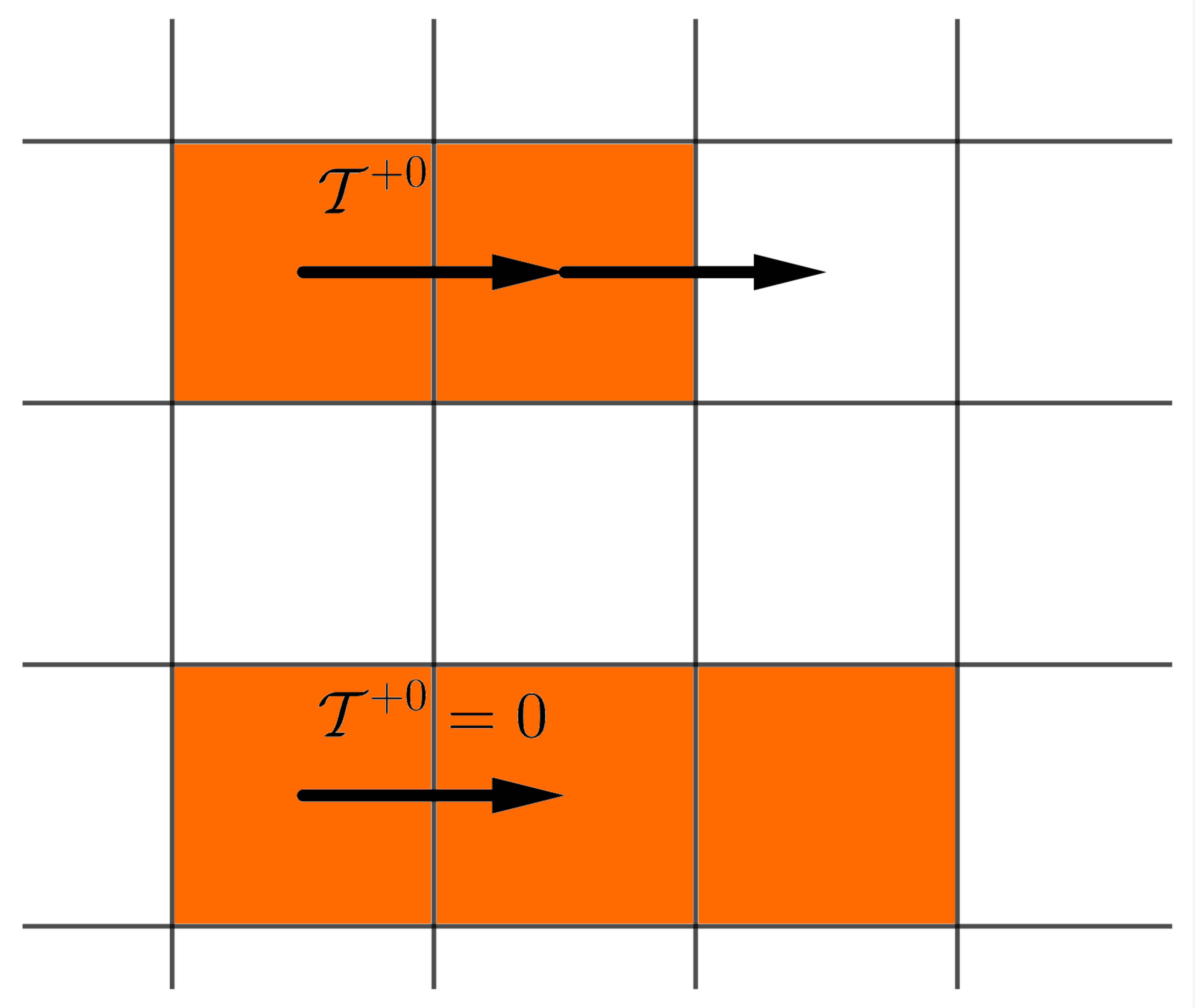}
\caption{Local pushing mechanism. Top: pushing is possible since the cell on the third column in free. Bottom: pushing is forbidden since the next site is occupied.}
\label{fig:pushing}
\end{subfigure}
\hfill
\begin{subfigure}[t]{0.45\textwidth}
\centering
\includegraphics[height=4cm]{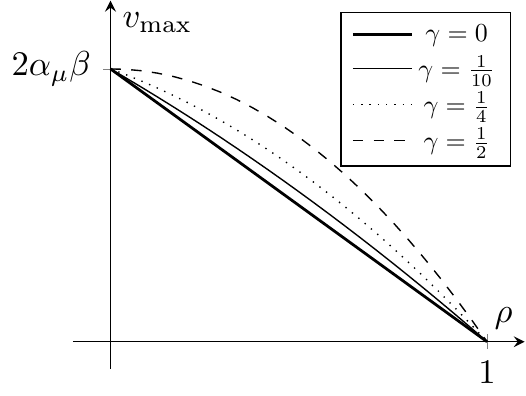}
\caption{Impact of pushing on the average velocity. Pushing leads to larger and increasingly concave velocities. }
\label{fig:fundamental}
\end{subfigure}
\caption{Effects of pushing}
\label{The effect of pushing}
\end{figure}

We will discuss the underlying CA approach for the sake of readability in 1D only, since its generalization to 2D is obvious. We assume that individuals
agents can move to a neighboring occupied cell with a given probability, by pushing the neighbor one cell further, provided that it is free.
Otherwise, such a move is forbidden. This mechanism is illustrated in Figure~\ref{fig:pushing}.
In 1D, the previously introduced transition rates are given by
\begin{align*}
\mathcal{T}^{i}(x)=\alpha_\mu \exp(\beta(\phi (x)-\phi (x+i\Delta x))).
\end{align*}
Since individuals can move to the right and left only, we will replace the superscript $i$ by $\pm$ indicating a jump to the respective neighboring sites.
Then the master equation in 1D is given by (ignoring constants):
\begin{equation*}
\begin{aligned}
    \rho(x,t+\Delta t)&-\rho(x,t)=\\[2mm] &-\rho (x)\mathcal{T}^{+}(x)\left( (1-\rho(x+\Delta x))+\gamma_\mu \rho(x+\Delta x) (1-\rho (x+2\Delta x)) \right)\\[2mm]
& -\rho (x)\mathcal{T}^{-}(x)\left( (1-\rho(x-\Delta x))+\gamma_\mu \rho(x-\Delta x) (1-\rho (x-2\Delta x)) \right)\\
\end{aligned}
\end{equation*}\begin{equation}
\begin{aligned}&+\rho(x+\Delta x)\mathcal{T}^{-}(x+\Delta x)(1-\rho(x))\\[2mm]
&+\gamma_\mu\rho(x+\Delta x)\rho(x+2\Delta x)\mathcal{T}^{-}(x+2\Delta x)(1-\rho(x))\\[2mm]
&+\rho(x-\Delta x)\mathcal{T}^{+}(x-\Delta x)(1-\rho(x))\\[2mm]
&+\gamma_\mu\rho(x-\Delta x)\rho(x-2\Delta x)\mathcal{T}^{+}(x-2\Delta x)(1-\rho(x))\,,
\end{aligned}
\label{eq:mepush}
\end{equation}
in which we omit $t$ for the sake of readability.
Here, $\gamma_\mu = \gamma(\mu) \in [0,1]$ denotes an increasing function and corresponds to the probability of an agent pushing. Note that we obtain
the original master equation~\eqref{eq:master} for $\gamma_\mu=0$.

\subsubsection*{Mean-field limit}
Using a formal Taylor expansion, we derive the limiting mean-field PDE where we generalized to 2D the approach mentioned previously:
\begin{align}\label{eq:pdepush}
\begin{aligned}
\partial_t \rho(x,y,t) &=\alpha_\mu \Div\bigl((1+4 \gamma_{\mu} \rho) \nabla \rho + 2 \beta \rho (1-\rho) (1+2\gamma_{\mu} \rho) \nabla \phi \bigr)\\[2mm]
\rho(0,x)&=\rho_0(x,y).
\end{aligned}
\end{align}
Equation \eqref{eq:pdepush} is supplemented with no-flux and outflow conditions of type \eqref{eq:bdc} for a modified flux
$\mathbf{j} = (1+4 \gamma_{\mu}\rho) \nabla \rho + 2 \beta \rho(1-\rho)(1+2\gamma_{\mu} \rho) \nabla \phi$ and $\alpha_{\mu}$ given by \eqref{eq:alpha}.

This equation has again a formal gradient flow structure with respect to the Wasserstein metric. The respective mobility and entropy are given by
\begin{align}\label{eq:pushmob}
m(\rho )=\alpha_\mu \rho (1-\rho) (1+2\gamma_{\mu} \rho),
\end{align}
and
\begin{multline}\label{eq:entropyfunctional2}
E(\rho ) = \int_{\Omega} \Big[ \frac{4 \gamma +1}{2 \gamma +1} ( 1-\rho) \log (1-\rho )\\[2mm]
 + \rho \log \rho +\frac{2 \gamma \rho +1}{2 \gamma +1} \log (2 \gamma \rho +1) +2 \beta \rho \phi\Big] ~ dx\,.
\end{multline}
We observe that the local pushing increases the mobility and the average velocity, see Figure \ref{fig:fundamental}. Furthermore, the velocity
decreases less in low density regimes and for higher motivation levels. Note that in case of pushing, the average velocity is always larger.

\begin{figure}[htb]
\vspace*{10pt} \centering
\begin{subfigure}[t]{0.8\textwidth}
\centering
 \includegraphics[width=0.8\textwidth]{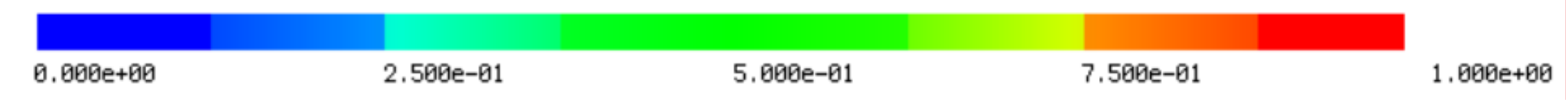}
\end{subfigure}

\begin{subfigure}[t]{0.4\textwidth}
\centering
 \includegraphics[height=4cm]{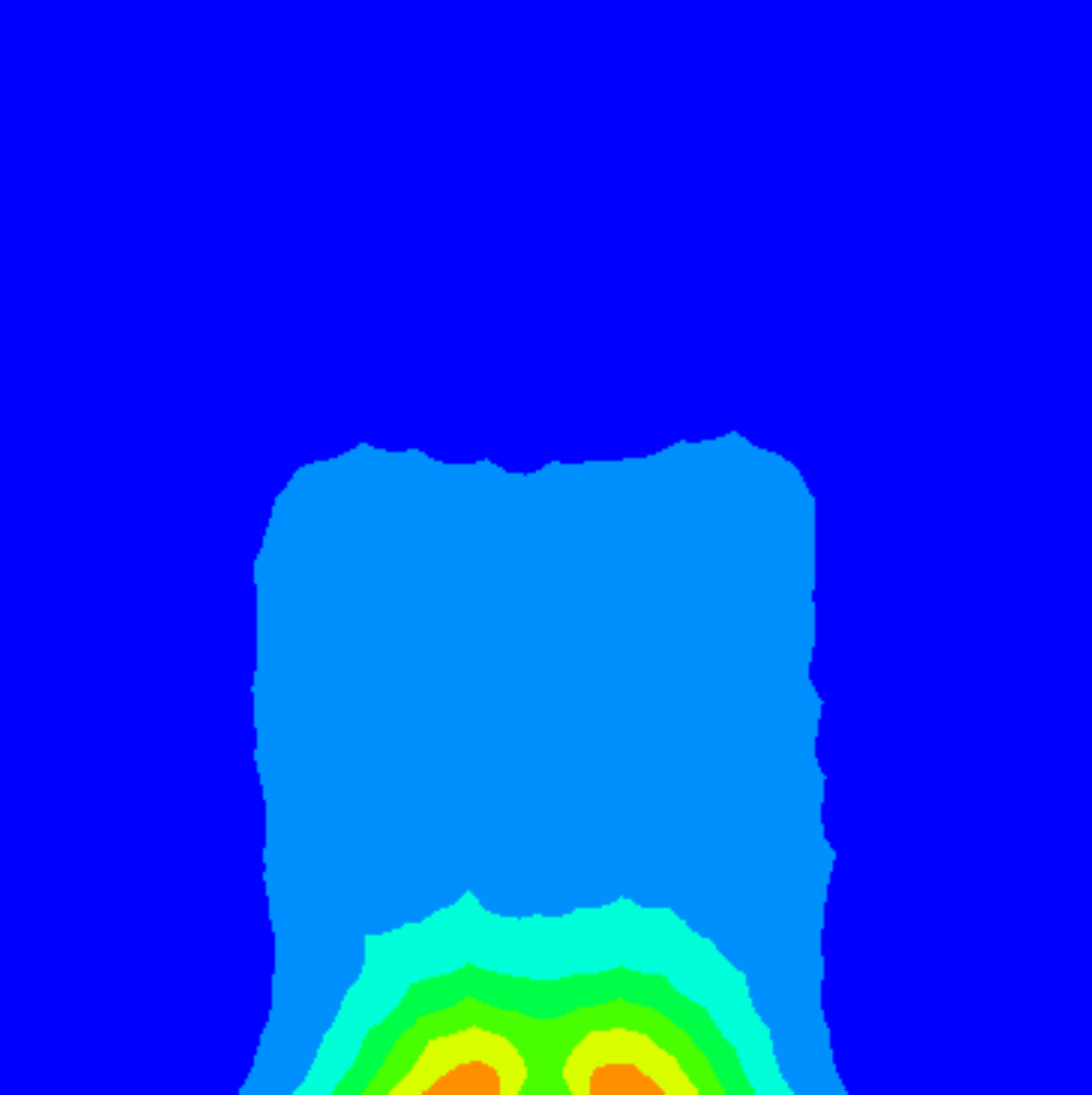}
 \caption{$t=5s$, $\gamma =0$, no pushing}
\end{subfigure}
\begin{subfigure}[t]{0.4\textwidth}
\centering
\includegraphics[height=4cm]{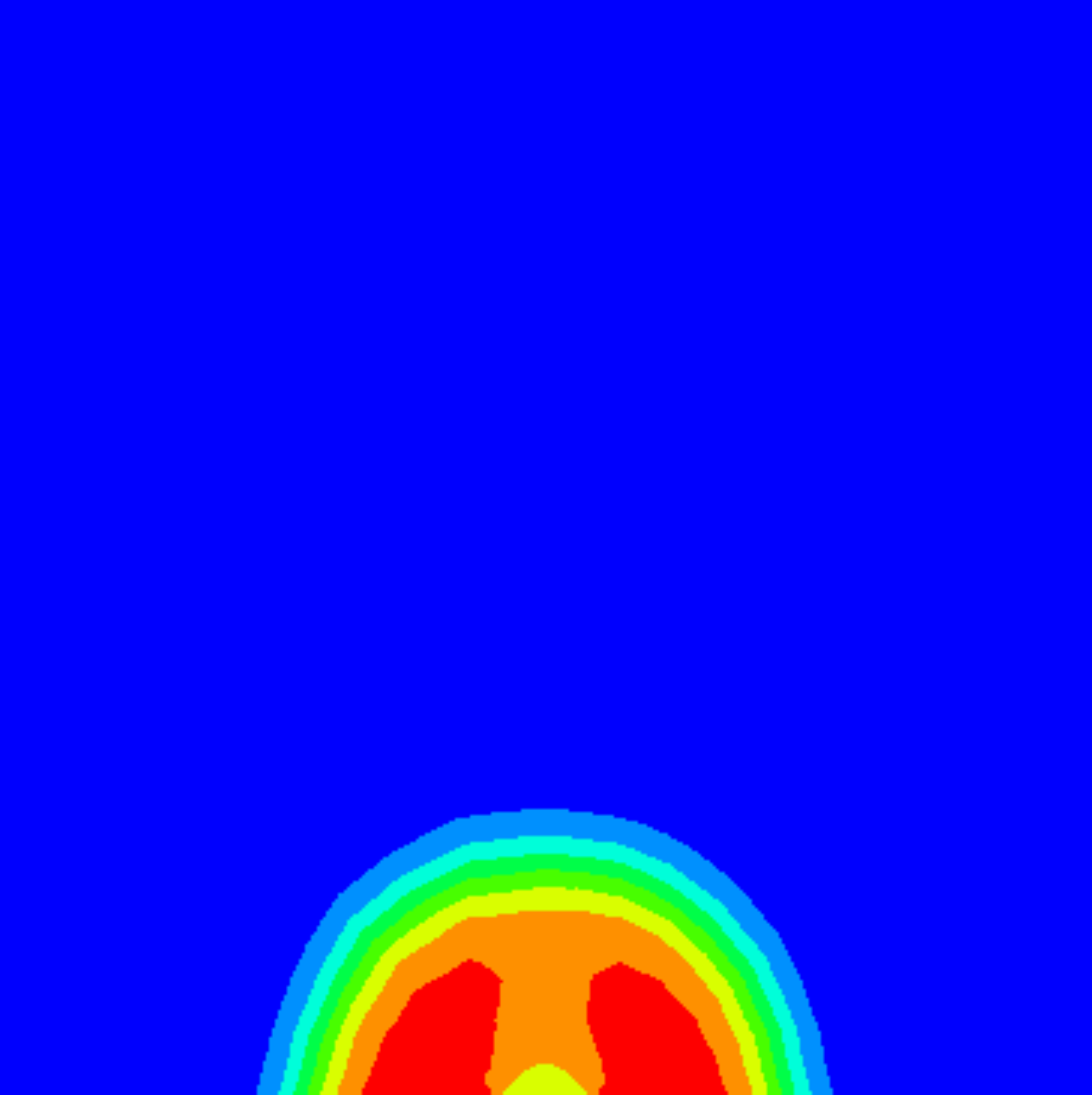}
\caption{$t=15s$, $\gamma =0$, no pushing}
\end{subfigure}

\begin{subfigure}[t]{0.4\textwidth}
\centering
 \includegraphics[height=4cm]{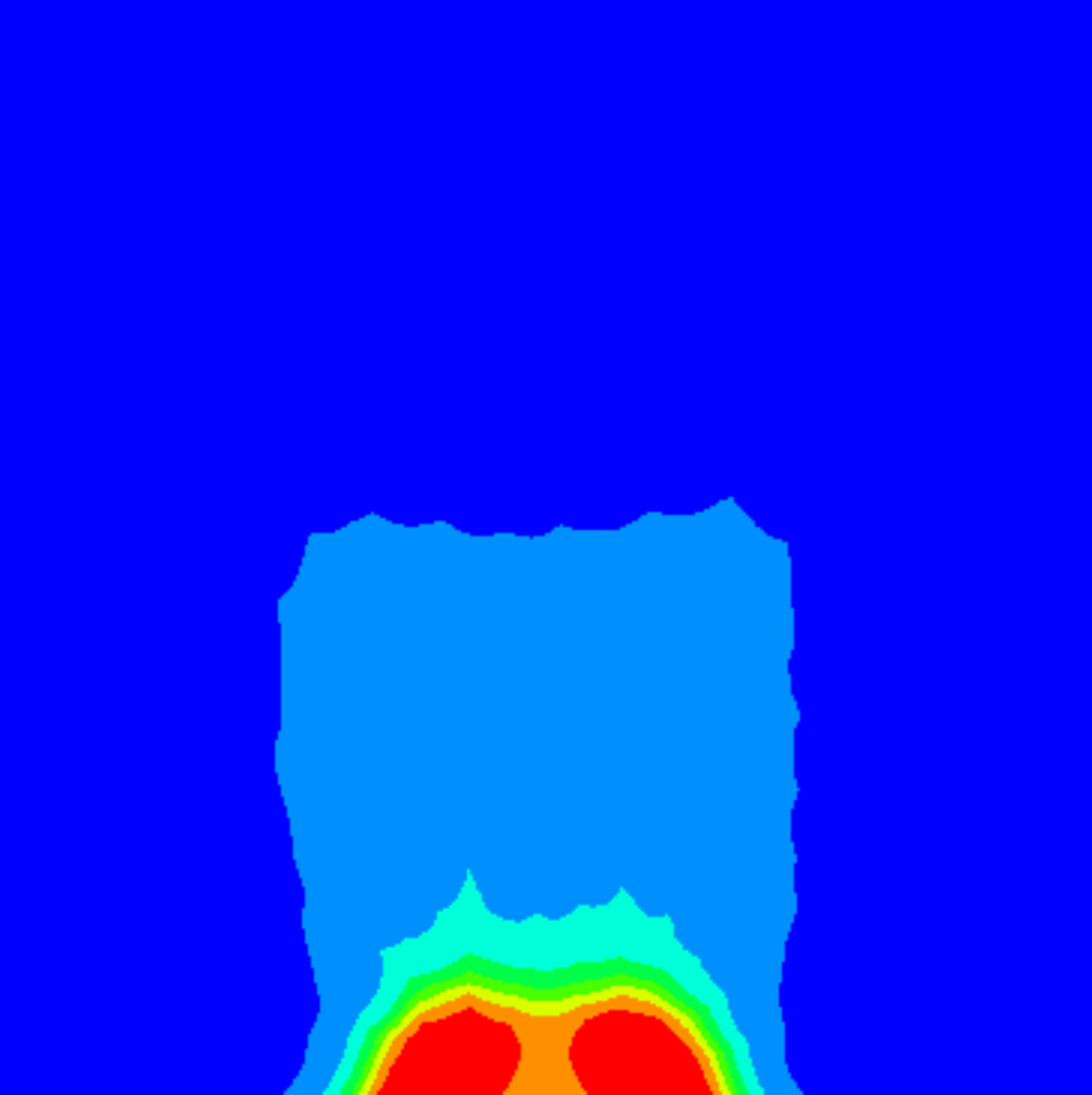}
 \caption{$t=5s$, $\gamma =1$, pushing}
\end{subfigure}
\begin{subfigure}[t]{0.4\textwidth}
\centering
\includegraphics[height=4cm]{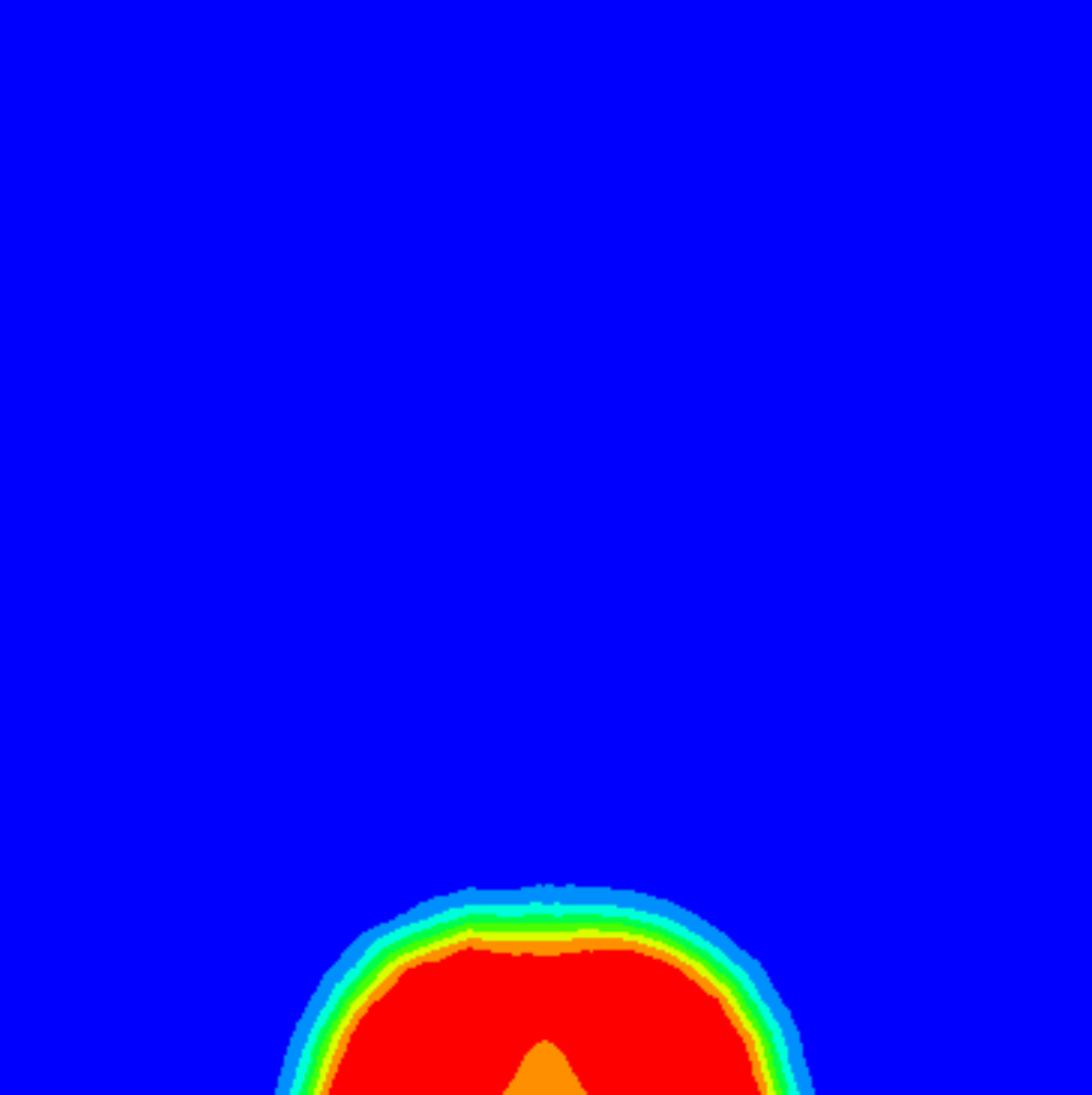}
\caption{$t=15s$, $\gamma =1$, pushing}
\end{subfigure}
\vspace*{10pt}\caption{Comparison of the congestion at the exit in case of pushing (bottom row)  and no-pushing (top row) for $60$ individuals:
We observe that people move faster towards the exit and the formation of a larger congested area in front of it.}\label{fig:comp_push}
\end{figure}

The local pushing weighs the $(1-\rho)\log(1-\rho)$ term and subsequently the finite volume effects much more. Furthermore,
it increases the entropy by an additional strictly positive term. Hence we expect a faster equilibration speed compared to the model of Section~\ref{sec:pde analysis}.
The expected behavior is confirmed by macroscopic simulations, see Figure \ref{fig:comp_push}. We consider a corridor filled with
$60$ people, and where we set $\beta = 3.8$. We observe that the individuals move faster towards the exit and that the congested area in front
of the exit builds up faster.

 Again, we recover the original PDE model by setting $\gamma_{\mu}=0$. The proof of global existence to~\eqref{eq:pdepush} follows arguments similar as for the
original PDE model \eqref{e:fpe}, since the pushing corresponds to a multiplicative prefactor in the mobility and a positive term in the entropy (which can be bounded).

\section{Conclusion}\label{s:con}
In this paper, we discussed micro- and macroscopic models for crowding and queuing at exits and bottlenecks, which were motivated by experiments conducted at the University
in Wuppertal. These experiments indicated that the geometry, ranging from corridors to open rooms, as well as the motivation level, such as a higher incentive to get to the exit
due to rewards, changes the overall dynamics significantly.

We propose a cellular automaton approach, in which the individual transition rates increase with the motivation level, and derive the corresponding continuum description using
a formal Taylor expansion. We use experimental data to calibrate the model and to understand the influence of parameters and geometry on the overall dynamics. Both the micro- and the macroscopic
description reproduce the experimental behavior correctly. In particular, we observe that corridors lead to lower densities and that the geometry has a stronger effect than the
motivation level. We plan to investigate the analysis of the coupled Hughes type models as well as the dynamics in case of pushing in more detail in the future.

%For acknowledgements section, please don't number the section, please begin it with \section*{Acknowledgements}
\section*{Acknowledgments}
The authors would like to thank Christoph Koutschan for the helpful discussions and input concerning the derivation of the respective mean field models using symbolic techniques.
Furthermore we would like to thank the team at the Forschungszentrum J\"ulich and the University of Wuppertal, in particular Armin Seyfried and Ben Hein, for providing the data and patiently answering all our questions.

All authors acknowledges partial support from the Austrian Academy of Sciences
via the New Frontier's grant NST 0001 and the EPSRC by the grant EP/P01240X/1.

\section*{Appendix}

\begin{proof}[Proof 1 of Theorem~\ref{theo:reg_eikonal}]
We start by a recalling a standard existence and regularity result from the literature, see~\cite{BCD2008} and ~\cite{a:cuc}.
Solutions to the eikonal equation~\eqref{eq:eik} in $\mathbb{R}^2\setminus \Gamma_E$ are given by the distance function
\begin{align*}
d(x,\Gamma_E)=\inf_{b\in\Gamma_E} |x-b|.
\end{align*}
Hence we discuss the regularity of $d$ in the following only.
We define the set
\begin{align*}
M(x)=\operatorname*{arg\,min}_{b\in\Gamma_E} d(x,\Gamma_E).
\end{align*}
If $\Gamma_E$ is a straight, bounded line, the $M$ is nonempty and consists of a single point for every $x\in\mathbb{R}^2$.
Since $|(\cdot-b)|$ is uniformly differentiable in $b$, and $b \mapsto D_x \lvert x - b \rvert$ is continuous,
we can deduce that the set $Y$
\begin{align*}
Y(x):=\{ D_x |x-b|: b\in M(x)\}
\end{align*}
is a singleton too. We now can apply Proposition 2.13 in~\cite{BCD2008} which states that $d$ is differentiable at $x$ if and only if $Y(x)$ is a singleton. Thus $d$ is differentiable for $\mathbb{R}^2\setminus \Gamma_E$, to be more precisely, we have $d\in C^1(\mathbb{R}^2\setminus \Gamma_E)\cap C(\mathbb R^2)$.

 Next we restrict $d$ to the corridor $\Omega\subset \mathbb R^2$ (being an open and bounded subset of $\Omega$).

Hence, $\phi_E \in C(\bar \Omega)\cap C^1 (\Omega)$. Since the $L^2$-norm of the first derivative of $\phi_E$ is
bounded by the equation (\ref{eq:eik}) itself, we can deduce that $\phi_E \in H^1(\Omega )$ since
\begin{align*}
\| \phi_E \|_{H^1(\Omega)}= \int_\Omega \phi_E^2dx+\int_\Omega (D\phi_E)^2dx \leq |\Omega | (\max \phi_E+1).
\end{align*}\qedhere
\end{proof}